%% file: paper.tex
\documentclass[10pt]{article}

\usepackage[hmargin=35mm,vmargin=40mm]{geometry}

\usepackage{amsmath}      
\usepackage{amssymb}      
\usepackage{amsthm}       
\usepackage{enumerate}    
\usepackage{graphicx}     
\usepackage{url}          

\newtheorem{theorem}{Theorem}[section]
\newtheorem{corollary}[theorem]{Corollary}
\newtheorem{lemma}[theorem]{Lemma}
\newtheorem{algorithm}[theorem]{Algorithm}

\newtheorem{defn}[theorem]{Definition}

\numberwithin{equation}{section}

\newcommand{\bdry}{\partial}

\newcommand{\euc}[1]{\mathbb{E}_{#1}}

\newcommand{\oproj}[1]{\mathbb{T}_{#1}}

\newcommand{\R}{\mathbb{R}}

\newcommand{\co}{\colon\thinspace}

\title{Projective geometry and the outer approximation \\
    algorithm for multiobjective linear programming}
\author{Benjamin A.~Burton and Melih Ozlen}

\date{June 15, 2010}

\begin{document}

\maketitle

\begin{abstract}
    A key problem in multiobjective linear programming is to find the
    set of all efficient extreme points in objective space.
    In this paper we introduce oriented projective geometry
    as an efficient and effective framework for solving this problem.
    The key advantage of oriented projective geometry is that we can
    work with an ``optimally simple'' but unbounded efficiency-equivalent
    polyhedron, yet apply to it the familiar theory and algorithms that
    are traditionally restricted to bounded polytopes.
    We apply these techniques to Benson's outer approximation algorithm,
    using oriented projective geometry to remove an exponentially large
    complexity from the algorithm and thereby remove a significant burden
    from the running time.


\end{abstract}

%
%

\input{intro.tex}
\input{prelim.tex}

\input{proj.tex}
\input{eff.tex}
\input{outer.tex}

\input{conc.tex}

%
%

\section*{Acknowledgements}

The first author is supported by the Australian Research Council
under the Discovery Projects funding scheme (project DP1094516).

%
%

\small
\bibliographystyle{amsplain}
\bibliography{pure}

%
%

\bigskip
\noindent
Benjamin A.~Burton \\
School of Mathematics and Physics, The University of Queensland \\
Brisbane QLD 4072, Australia \\
(bab@debian.org)

\bigskip
\noindent
Melih Ozlen \\
School of Mathematical and Geospatial Sciences, RMIT University \\
GPO Box 2476V, Melbourne VIC 3001, Australia \\
(melih.ozlen@rmit.edu.au)

\end{document}

%% file: intro.tex
\section{Introduction} \label{s-intro}

In traditional linear programming, our goal is to maximise a linear
objective function over some polytope $X \subseteq \euc{n}$ (where
$\euc{n}$ denotes $n$-dimensional Euclidean space).
In \emph{multiobjective linear programming}, we attempt to maximise many
linear objective functions over $X$ simultaneously.  Multiobjective linear
programming is an important practical tool in multicriteria decision making
\cite{cohon78-molp,stadler88-molp,steuer86-molp,zeleny82-molp},
and recent decades have seen a surge of activity in understanding the
complexities of such problems and developing algorithmic tools to solve
them.

Typically different objectives achieve their global maxima at different
points of the polytope $X$, resulting in trade-off situations in which
some objectives can be improved at the expense of others.
In the most general solution to a multiobjective linear program, we seek to
find all \emph{efficient} (or \emph{Pareto-optimal}) points, which are those
solutions from which we cannot increase any one objective without
simultaneously decreasing another \cite{geoffrion68-efficiency}.

The set of all efficient points in $X$ is difficult to manage, since it
is typically large and non-convex.  Instead it is preferable to seek a
smaller combinatorial description of the efficient set.  A common
approach (which we follow in this paper) is to compute the
\emph{efficient extreme points} of $X$ (that is, the efficient vertices
of $X$), from which the full efficient set can be reconstructed
\cite{ecker80-faces,gal77-general,isermann77-enumeration}.

Early approaches to multiobjective linear programming were based on
techniques from the simplex method
\cite{armand91-determination,evans84-overview,yu75-nondominated,zeleny82-molp},
and involved walking through the polytope $X$.
However, for problems with large numbers of variables $n$, the polytope
$X \subseteq \euc{n}$ can grow to become extremely complex.

A new approach was therefore developed to work in the smaller
\emph{objective space} $\euc{p}$, where $p$ is the number of objectives
under consideration
\cite{dauer90-constructing,dauer87-objective,dauer90-solving}.
The original polytope $X \subseteq \euc{n}$ naturally maps to an
\emph{outcome set} $Y^= \subseteq \euc{p}$, which is a polytope
that describes all achievable combinations of objectives.

Working in objective space is preferable because the running time
of many relevant algorithms is exponential in the dimension of
the polytope \cite{steuer05-regression,ziegler95}---this means that
working in $p$ dimensions instead of $n$
can have an enormous impact on the computational resources required.
Moreover, it is often easier for the decision maker to choose a preferred
outcome from $Y^=$ instead of a specific combination of input variables
from $X$ \cite{benson98-outer,dauer90-constructing}.

We therefore focus on the problem of finding all efficient extreme
points of the outcome set $Y^=$.
There are many algorithms in the literature for solving this problem.
One notable family of algorithms is based on weight space decompositions
\cite{benson00-partition,benson02-decomposition,przybylski09-recursive},
extending earlier theoretical work on the vector maximum problem
\cite{isermann74-efficiency,yu75-nondominated}.
Another important family of algorithms is based on outer approximation
techniques from the global optimisation literature
\cite{benson98-analysis,benson98-hybrid,benson98-outer}.
It remains relatively unknown how these algorithms
compare in practice, and it is possible that both families have important
roles to play as the problem size becomes large \cite{benson02-decomposition}.

A key optimisation in these algorithms comes from working with
\emph{efficiency-equivalent polyhedra}.  These are polytopes and
polyhedra in $\euc{p}$ whose efficient points are the same as for
the outcome set $Y^=$, but whose non-efficient regions are much simpler.
Two such polyhedra feature prominently in the literature:
\begin{itemize}
    \item Dauer and Saleh \cite{dauer90-constructing,dauer92-representation}
    define the efficiency-equivalent polyhedron $Y^\leq$,
    which is combinatorially simple but which has the practical
    disadvantage of being unbounded
    (extending infinitely far in the negative directions).
    Being unbounded means that it cannot be expressed as the convex hull
    of its vertices, which causes problems for important procedures such
    as the double description method of Motzkin et~al.\ \cite{motzkin53-dd}.

    \item Benson \cite{benson98-analysis,benson98-hybrid,benson98-outer}
    describes the efficiency-equivalent polytope $Y^\square$,
    which is bounded with a cube-like structure.
    Although it is less combinatorially simple,
    the bounded nature of $Y^\square$ makes it more practical for use in
    algorithms, and indeed $Y^\square$ plays a key role in the family of outer
    approximation algorithms mentioned above.
\end{itemize}

In this paper we develop a new framework for algorithms in
multiobjective linear programming, based on \emph{oriented projective
geometry}.  Introduced by Stolfi \cite{stolfi87-oriented}, oriented
projective geometry augments Euclidean space with additional points at
infinity.  Unlike classical projective geometry, oriented projective
geometry also preserves notions such as convexity, half-spaces and
polytopes, all of which are crucial for linear programming.

We show that by working in oriented projective $p$-space, we can achieve
the best of both worlds with efficiency-equivalent polyhedra.
Specifically, we work with an efficiency-equivalent \emph{projective
polytope} $T^\leq$, which inherits both the simple combinatorial
structure of $Y^\leq$ and the wider practical usability of $Y^\square$.

To show the \emph{efficiency} of this new framework, we undertake a
combinatorial analysis of the bounded Euclidean polytope $Y^\square$
and compare this to the projective polytope $T^\leq$.
We show that the non-efficient structure (that is, the ``uninteresting''
region) of $Y^\square$ is exponentially complex: the number of
non-efficient vertices is always at least $2^p-1$, and can grow
as fast as $\Omega(n^{\lfloor (p-1)/2 \rfloor})$.
In contrast, we show that the projective polytope $T^\leq$ has just
$p$ non-efficient vertices, a figure that is not only small but also
independent of the problem size $n$.

Even for problems where the number of objectives $p$ is small, our
framework enjoys a significant advantage by using a projective
polytope whose non-efficient region has a complexity independent of $n$.
For problems where $p$ grows larger,
such as the aircraft control problem of Schy and Giesy with 70 objectives
\cite{schy88-aircraft}, or applications in welfare economics where
$p$ can grow arbitrarily large \cite{schulz88-welfare}, even the lower
bound of $2^p-1$ puts $Y^\square$ at a severe disadvantage compared to
the projective polytope $T^\leq$ that we use in our framework.

To show the \emph{effectiveness} of this framework, we apply it
to the outer approximation algorithm of Benson \cite{benson98-outer}.
The result is a new algorithm that works in oriented projective space,
combining the practical advantages of the original algorithm
with the combinatorial simplicity of $T^\leq$.
By studying the inner workings of this algorithm, we show how the
complexity results above translate into significant benefits in terms
of running time.

The layout of this paper is as follows.
Section~\ref{s-prelim} summarises relevant
definitions and results from polytope theory and multiobjective linear
programming, and Section~\ref{s-proj} offers a detailed introduction to
oriented projective geometry.  In Section~\ref{s-eff} we discuss
efficiency-equivalent polyhedra, and prove the combinatorial complexity
results outlined above.
We apply our techniques to the outer approximation algorithm
in Section~\ref{s-outer}, and in Section~\ref{s-conc} we conclude with a
discussion of the theoretical and practical implications of our work.

%% file: prelim.tex
\section{Preliminaries} \label{s-prelim}

In this brief section, we outline key definitions from the theory of
polytopes and polyhedra and the theory of multiobjective linear
programming.  For polytopes and polyhedra we follow the terminology of
Ziegler \cite{ziegler95} and
Matou{\v s}ek and G{\"a}rtner \cite{matousek07-linear}.
For multiobjective linear programming we follow the conventions and notation
used by Benson \cite{benson98-outer}.

Throughout this paper we work in both traditional Euclidean geometry and
the alternative \emph{oriented projective geometry}.
We denote $d$-dimensional
Euclidean space by $\euc{d}$, and each point in $\euc{d}$
is described in the usual way by a coordinate vector
$\mathbf{x} = (x_1,\ldots,x_d) \in \R^d$.
Although the geometry $\euc{d}$ and the coordinate space $\R^d$ are
identical in the Euclidean case, this distinction becomes important in
oriented projective geometry, where we describe $d$-dimensional geometry
using coordinates in $\R^{d+1}$.  We delay any further discussion
of oriented projective geometry until Section~\ref{s-proj}.

The concepts of polytopes and polyhedra are fundamental to linear programming.
A \emph{polyhedron} in $\euc{d}$ is an intersection of finitely many
closed half-spaces in $\euc{d}$; this may be either bounded (like a
cube) or unbounded (like an infinite cone).
A polyhedron that is bounded is also called a \emph{polytope}.
Every polytope in $\euc{d}$ can be expressed as a convex hull of
finitely many vertices in $\euc{d}$, whereas unbounded polyhedra cannot
be expressed in this way.

A \emph{hyperplane} in $\euc{d}$ is a $(d-1)$-dimensional affine subspace,
and consists of all points with coordinates
$\{\mathbf{x} \in \R^d\,|\,\mathbf{h}\cdot\mathbf{x}=\alpha\}$
for some $\mathbf{h} \in \R^{d}$ and $\alpha \in \R$.
If $P \subseteq \euc{d}$ is a polyhedron (bounded or unbounded)
and $H \subseteq \euc{d}$ is a hyperplane, then we call
$H$ a \emph{supporting hyperplane} for $P$ if
(i)~$H$ contains some point of $P$, and (ii)~all of the points in $P$
lie within $H$ and/or to the \emph{same side} of $H$.
In other words,
(i)~$\mathbf{h} \cdot \mathbf{x} = \alpha$ for some $\mathbf{x} \in P$,
and
(ii)~either $\mathbf{h} \cdot \mathbf{x} \geq \alpha$ for all
$\mathbf{x} \in P$, or else
$\mathbf{h} \cdot \mathbf{x} \leq \alpha$ for all $\mathbf{x} \in P$.
The intersection of $P$ with a supporting hyperplane is called a
\emph{face} of $P$.\footnote{For theoretical convenience, two additional
    faces are defined: a ``full face'' containing all of $P$, and an
    ``empty face'' containing no points at all.  However, neither of
    these is relevant to this paper.}

The \emph{dimension} of a polyhedron is defined to be the dimension of its
affine hull (in particular, a polyhedron in $\euc{d}$ must have dimension
$\leq d$).  Likewise, the dimension of a face is defined to be the
dimension of its affine hull.  If $P$ is a polyhedron of dimension $k$,
then faces of dimension $0$, $1$ and $(k-1)$ are called
\emph{vertices}, \emph{edges} and \emph{facets} of $P$ respectively.

A \emph{multiobjective linear program} requires us to maximise
$p$ linear objective functions over a polyhedron in $\euc{n}$.
This polyhedron is called the \emph{feasible solution set} $X$,
and consists of all $\mathbf{x} \in \euc{n}$ for which
$A \mathbf{x} = \mathbf{b}$ and $\mathbf{x} \geq \mathbf{0}$,
where $A$ is an $m \times n$ matrix and $\mathbf{b} \in \R^m$.
Like Benson \cite{benson98-outer}, we assume for convenience that
$X$ is bounded; that is, $X$ is a polytope.

The $p$ linear objectives are described by a $p \times n$
\emph{objective matrix} $C$, and our task is to maximise
$C \mathbf{x}$ over all $\mathbf{x} \in X$.
We define the \emph{outcome set} $Y^=$ as
$\{ C\mathbf{x}\,|\,\mathbf{x} \in X\}$.
The outcome set $Y^=$ is in turn a polytope in $\euc{p}$, and our task
can be restated as maximising $\mathbf{y} \in Y^=$.

In general we cannot achieve a global maximum for all $p$ objectives
simultaneously, so instead we focus on non-dominated outcomes.
We call an outcome $\mathbf{y} \in Y^=$ \emph{efficient}
(also called \emph{Pareto-optimal} or \emph{non-dominated})
if there is no other $\mathbf{y}' \in Y^=$
for which $\mathbf{y}' \geq \mathbf{y}$.
The \emph{efficient outcome set}
$Y^=_E$ is defined to be the set of all efficient outcomes in $Y^=$.

More generally, an efficient point in any polytope or polyhedron
$Z \subseteq \euc{p}$ is a point $\mathbf{z} \in Z$ for which there is
no other $\mathbf{z}' \in Z$ satisfying $\mathbf{z}' \geq \mathbf{z}$.
The set of all efficient points in $Z$ is denoted $Z_E$,
and this set can always be expressed as a union of faces of $Z$.
If $Z_E = Y^=_E$ then we call $Z$ an \emph{efficiency-equivalent
polyhedron} for $Y^=$.

Since the efficient outcome set $Y^=_E$ is generally infinite,
we focus our attention on the set of efficient extreme outcomes.
An \emph{extreme point} of a polyhedron is simply a vertex of the
polyhedron, and an \emph{efficient extreme outcome} is an
efficient vertex of $Y^=$.  Because $Y^=$ is a polyhedron there are only
finitely many efficient extreme outcomes, and from these we can generate
all of $Y^=_E$ if we so desire \cite{gallagher95-representation}.
Our focus then for the remainder of this paper is on generating the
finite set of efficient extreme outcomes.

%% file: proj.tex
\section{Oriented projective geometry} \label{s-proj}

In \emph{classical projective geometry}, we augment the Euclidean space
$\euc{d}$ with additional ``points at infinity''.  This offers
significant advantages over traditional
Euclidean geometry.  Theoretically, we achieve a powerful
duality between points and hyperplanes.
Computationally, we can perform simple arithmetic on infinite
limits (such as the ``endpoints'' of lines and rays); moreover, we
lose the notion of ``parallelness'' and thereby avoid a myriad of
special cases.
For a thorough overview of classical projective geometry, the reader is
referred to texts such as Coxeter \cite{coxeter87-projective} or
Beutelspacher and Rosenbaum \cite{beutelspacher98-projective}.

A key drawback of classical projective geometry is that we
also lose important concepts such as convexity, half-spaces and polytopes.
\emph{Oriented projective geometry} further augments classical
projective geometry to restore these concepts.
For our purposes, the primary advantage of
oriented projective geometry is that it allows us to treat an
\emph{unbounded} polyhedron just like a bounded polytope---that is,
as a convex hull of finitely many vertices (where some of these vertices
happen to be points at infinity).

Oriented projective geometry was introduced by Stolfi in the 1980s
\cite{stolfi87-oriented,stolfi91-oriented}, and has since found
applications in computer vision \cite{laveau96-oriented,werner01-matching}.
Kirby presents a leisurely geometric overview in \cite{kirby02-celestial},
and Boissonnat and Yvinec discuss issues relating to
polytopes and polyhedra in \cite{boissonnat98-geometry}.
In this section we describe the way in which oriented projective geometry
extends Euclidean geometry, show how we can perform arithmetic in oriented
projective geometry using \emph{signed homogeneous coordinates},
and discuss how unbounded polyhedra in Euclidean space can be
reinterpreted as \emph{projective polytopes}.

Following the original notation of Stolfi, we denote
the \emph{oriented projective $d$-space} by $\oproj{d}$.
Geometrically, we construct $\oproj{d}$ by augmenting the Euclidean
$d$-space $\euc{d}$ as follows:
\begin{itemize}
    \item We add a collection of \emph{points at infinity}.
    Each ray of $\euc{d}$
    passes through one such point, and each line of $\euc{d}$ passes through
    two such points (one in each direction).  Two rays pass through the
    same point at infinity if and only if they are parallel and
    pointing in the same direction.  This is illustrated for $d=2$ in
    Figure~\ref{fig-celestial}.

    \begin{figure}[htb]
        \centering
        \includegraphics[scale=0.9]{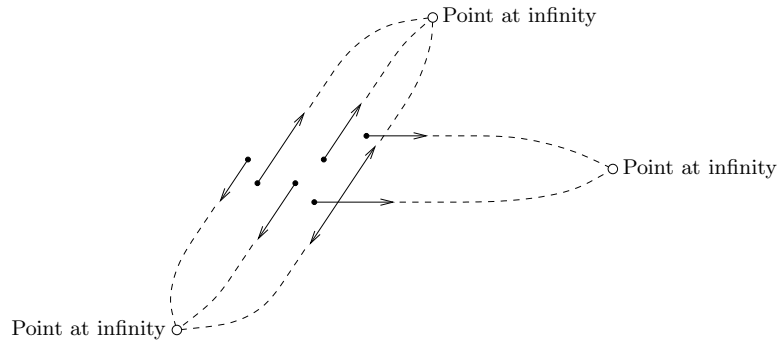}
        \caption{Euclidean rays and lines passing through points at
            infinity in $\oproj{2}$}
        \label{fig-celestial}
    \end{figure}

    \item Travelling beyond these points at infinity, we add a second
    ``hidden'' copy of Euclidean $d$-space which we call \emph{invisible space};
    in contrast, we refer to our original copy of $\euc{d}$ as
    \emph{visible space}.
    For computational purposes we can ignore invisible space---all
    of the Euclidean structures that we use in this
    paper will be situated in the visible copy of $\euc{d}$.
    Invisible space is simply a theoretical convenience that
    endows $\oproj{d}$ with a number of
    useful geometric and analytical properties.
\end{itemize}

\begin{figure}[htb]
    \centering
    \includegraphics[scale=0.9]{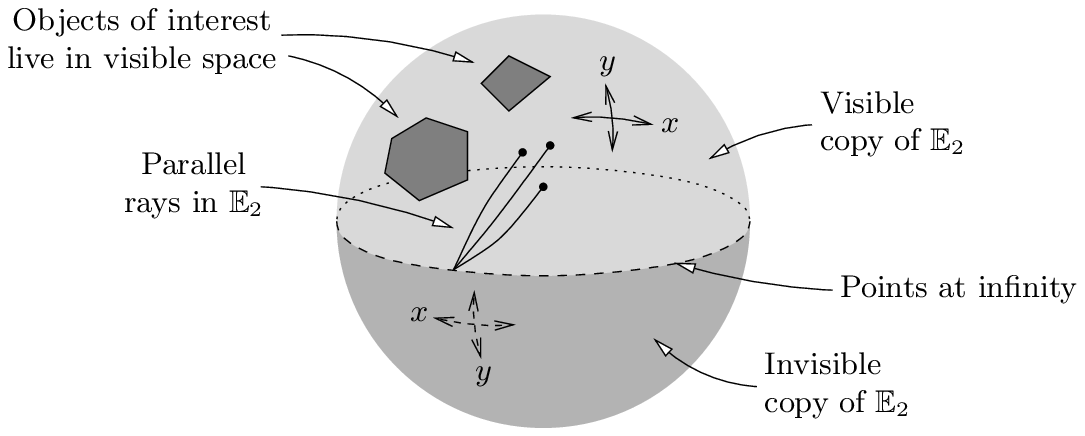}
    \caption{Building the oriented projective space $\oproj{2}$}
    \label{fig-proj-plane}
\end{figure}

For illustration, Figure~\ref{fig-proj-plane} shows the full
construction of the oriented projective space $\oproj{2}$.
We see from this diagram that $\oproj{2}$ has the
topological structure of a sphere: the visible points lie on the
northern hemisphere, the invisible points lie on the southern
hemisphere, and the points at infinity lie on the equator.
This picture generalises to higher dimensions, and
Stolfi's original construction shows how points in $\oproj{d}$
correspond to points on the unit sphere in $\R^{d+1}$ via a
geometric projection.

We coordinatise oriented projective $d$-space using
\emph{signed homogeneous coordinates}.  Each point $\mathbf{p}\in\oproj{d}$
is represented by a $(d+1)$-tuple of numbers, and two $(d+1)$-tuples
$\mathbf{x},\mathbf{y} \in \R^{d+1}$ represent the \emph{same}
point of $\oproj{d}$ if and only if $\mathbf{x} = \lambda \mathbf{y}$
for some scalar $\lambda > 0$.  So, for instance, the coordinates
$(2,3,1)$, $(4,6,2)$ and $(1,\frac32,\frac12)$ all represent the same
point of $\oproj{2}$, but the coordinates $(-2,-3,-1)$ represent a
different point.
The coordinates $(0,\ldots,0)$ do not represent any point of $\oproj{d}$
at all.

\begin{figure}[htb]
    \centering
    \includegraphics[scale=0.9]{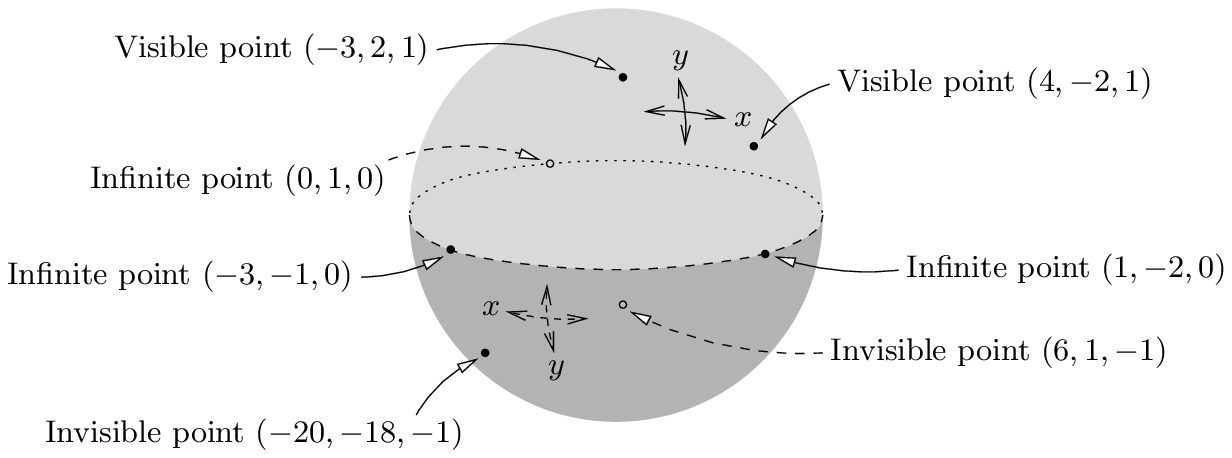}
    \caption{Different types of points in $\oproj{2}$}
    \label{fig-proj-points}
\end{figure}

More specifically,
we coordinatise the individual points of $\oproj{d}$ as follows:
\begin{itemize}
    \item Consider any Euclidean point
    $\mathbf{e} = (e_1,\ldots,e_d) \in \euc{d}$.
    In $\oproj{d}$, the copy of $\mathbf{e}$ in visible space is
    described by the homogeneous coordinates
    $(\lambda e_1,\ldots,\lambda e_d,\lambda)$ for any $\lambda > 0$,
    and the copy of $\mathbf{e}$ in invisible
    space is described by the homogeneous coordinates
    $(-\lambda e_1,\ldots,-\lambda e_d,-\lambda)$ for any $\lambda > 0$.

    \item Consider any Euclidean ray
    $\{\mu \mathbf{r}\,|\,\mu \geq 0\} \subseteq \euc{d}$, where
    $\mathbf{r}=(r_1,\ldots,r_d)$.
    In $\oproj{d}$, this ray in visible space extends to a point at
    infinity, which is described by the homogeneous coordinates
    $(r_1,\ldots,r_d,0)$.
\end{itemize}
Conversely, any point $\mathbf{x}\in\oproj{d}$ with homogeneous
coordinates $(x_1,\ldots,x_{d+1})$ can be interpreted geometrically as
follows:
\begin{itemize}
    \item If $x_{d+1}>0$ then $\mathbf{x}$ describes the visible
    Euclidean point
    $(\frac{x_1}{x_{d+1}},\ldots,\frac{x_d}{x_{d+1}}) \in \euc{d}$.

    \item If $x_{d+1}<0$ then $\mathbf{x}$ describes the invisible
    Euclidean point
    $(\frac{x_1}{x_{d+1}},\ldots,\frac{x_d}{x_{d+1}}) \in \euc{d}$.

    \item If $x_{d+1}=0$ then $\mathbf{x}$ describes a point at
    infinity, which is the limiting point of the visible Euclidean ray
    $\{\mu \cdot (x_1,\ldots,x_d)\,|\,\mu \geq 0\} \subseteq \euc{d}$.
\end{itemize}
For brevity, we refer to such points as \emph{visible},
\emph{invisible} and \emph{infinite} points of $\oproj{d}$ respectively.
Examples of all three types of point in $\oproj{2}$ are illustrated in
Figure~\ref{fig-proj-points}.

\begin{figure}[htb]
    \centering
    \includegraphics[scale=0.9]{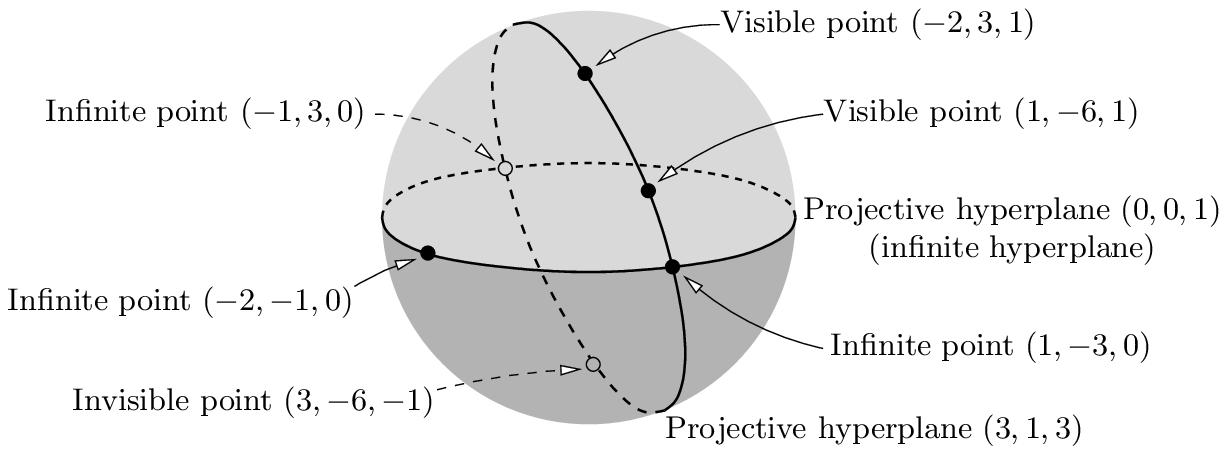}
    \caption{Two examples of projective hyperplanes in $\oproj{2}$}
    \label{fig-proj-hyp}
\end{figure}

A \emph{projective hyperplane} in $\oproj{d}$ is coordinatised by a
$(d+1)$-tuple $\mathbf{h}=(h_1,\ldots,h_{d+1}) \in \R^{d+1}$,
and contains all points
$\mathbf{x}\in\oproj{d}$ for which $h_1x_1+\ldots+h_{d+1}x_{d+1}=0$,
or more concisely $\mathbf{h}\cdot\mathbf{x}=0$.
Most projective hyperplanes
consist of a Euclidean hyperplane in visible space and a
Euclidean hyperplane in invisible space, joined together by a
common set of points at infinity.  The only exception
is the \emph{infinite hyperplane} described by the $(d+1)$-tuple
$(0,\ldots,0,1)$, which contains all of the points at infinity and
nothing else.  Projective hyperplanes of both types are illustrated in
Figure~\ref{fig-proj-hyp}.  The $(d+1)$-tuple $\lambda\mathbf{h}$
describes the same projective hyperplane as $\mathbf{h}$ for any
$\lambda \neq 0$, and the $(d+1)$-tuple $(0,\ldots,0)$ does not describe
any hyperplane at all.

\begin{figure}[htb]
    \centering
    \includegraphics[scale=0.9]{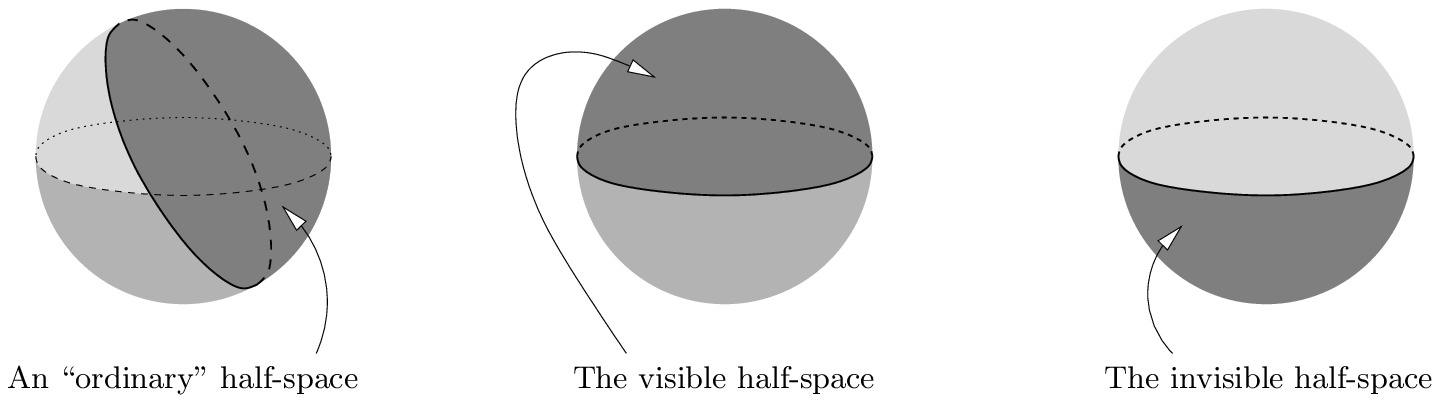}
    \caption{Different types of projective half-space in $\oproj{2}$}
    \label{fig-proj-half}
\end{figure}

A \emph{projective half-space} in $\oproj{d}$ is likewise coordinatised
by a $(d+1)$-tuple $\mathbf{h}=(h_1,\ldots,h_{d+1}) \allowbreak \in \R^{d+1}$,
and contains all points
$\mathbf{x}\in\oproj{d}$ for which $\mathbf{h}\cdot\mathbf{x} \geq 0$.
Most projective half-spaces consist of a Euclidean half-space in visible
space and a Euclidean half-space in invisible space, again joined together
by a common set of points at infinity.  There are two exceptions:
the \emph{visible half-space} described by the $(d+1)$-tuple
$(0,\ldots,0,1)$ and the \emph{invisible half-space} described by the
$(d+1)$-tuple $(0,\ldots,0,-1)$.
The visible half-space contains all visible and infinite points of
$\oproj{d}$, and the invisible half-space contains all invisible and
infinite points of $\oproj{d}$.
All three types of projective
half-space are illustrated in Figure~\ref{fig-proj-half}.
The $(d+1)$-tuple $\lambda\mathbf{h}$ describes the same projective
half-space as $\mathbf{h}$ for any
$\lambda > 0$ (but not $\lambda < 0$),
and as usual $(0,\ldots,0)$ does not describe any half-space at all.

For a projective half-space with coordinates
$\mathbf{h}=(h_1,\ldots,h_{d+1})$, the \emph{interior} of this
half-space consists of all points
$\mathbf{x}\in\oproj{d}$ for which $\mathbf{h}\cdot\mathbf{x} > 0$;
that is, all points of the half-space that do not lie on the boundary
hyperplane.  It is worth noting that the interior of the visible
half-space is precisely the set of all visible Euclidean points.

For any point $\mathbf{x} \in \oproj{d}$ with coordinates
$(x_1,\ldots,x_{d+1})$, the \emph{opposite} point to $\mathbf{x}$ is the
point with coordinates $(-x_1,\ldots,-x_{d+1})$, which we denote
$\neg\mathbf{x}$.  The opposite of a visible point is an invisible
point, and vice versa.  If $\mathbf{x}$ is a point at infinity then
$\neg\mathbf{x}$ is another point at infinity, and these points lie at
opposite ends of the same Euclidean lines.

Arithmetic in $\oproj{d}$ is very simple.  For distinct and
non-opposite points $\mathbf{x},\mathbf{y} \in \oproj{d}$ with coordinates
$(x_1,\ldots,x_{d+1})$ and $(y_1,\ldots,y_{d+1})$ respectively,
the unique \emph{line} joining $\mathbf{x}$ and $\mathbf{y}$ consists of all
points with coordinates
$(\alpha x_1 + \beta y_1,\ldots,\alpha x_{d+1} + \beta y_{d+1})$
for any $\alpha,\beta \in \R$, excluding the case $\alpha=\beta=0$.\footnote{%
    We insist on non-opposite points because there are infinitely many
    distinct lines in $\oproj{d}$ joining $\mathbf{x}$ and $\neg\mathbf{x}$.}
For the \emph{line segment} joining $\mathbf{x}$ and $\mathbf{y}$ we add
the constraints $\alpha,\beta \geq 0$.
In other words, the arithmetic of lines in $\oproj{d}$ is essentially the
same as traditional Euclidean arithmetic in $\R^{d+1}$, except that we
discard the constraint $\alpha+\beta=1$ (because signed homogeneous
coordinates are unaffected by positive scaling).

We can define convexity in terms of line segments.
Consider any set $S \subseteq \oproj{d}$ that lies within the interior
of some projective half-space.  Then $S$ is \emph{convex} if and only if,
for any two points $\mathbf{x},\mathbf{y} \in S$, the entire line segment
joining $\mathbf{x}$ and $\mathbf{y}$ lies within $S$.
The half-space constraint ensures that $S$ cannot contain
two opposite points, which means that line segments are always well-defined.
Note that this half-space constraint is automatically satisfied for
\emph{any} set of visible Euclidean points (which all lie in the
interior of the visible half-space).

Consider now some finite set of points
$\mathbf{x}^{(1)},\ldots,\mathbf{x}^{(k)} \in \oproj{d}$,
where each $\mathbf{x}^{(i)}$ has coordinates
$(x^{(i)}_1,\ldots,x^{(i)}_{d+1})$,
and again suppose that these points all lie within the interior of
some projective half-space.  The \emph{convex hull} of these points is
the smallest convex set containing all of
$\mathbf{x}^{(1)},\ldots,\mathbf{x}^{(k)}$.  Equivalently, it is
the set of all points in $\oproj{d}$ with coordinates of the form
$(\sum \alpha_i x^{(i)}_1,\ldots,\sum \alpha_i x^{(i)}_{d+1})$
for scalars $\alpha_1,\ldots,\alpha_i \geq 0$, excluding
the case $\alpha_1=\ldots=\alpha_i=0$.

It is important to note that, when we restrict our attention to
visible Euclidean points, all of these concepts---lines, line segments,
convexity and convex hulls---have precisely the same meaning as in
traditional Euclidean geometry.

As noted earlier, the advantage of oriented projective geometry
for us is the way in which it handles polytopes and polyhedra.
Recall that in Euclidean geometry, a polyhedron (which may be
bounded or unbounded) can be expressed as an intersection of finitely
many half-spaces, and that a polytope (that is, a bounded polyhedron)
can also be expressed as a convex hull of finitely many vertices.
Using oriented projective geometry, we are able to express
some \emph{unbounded} polyhedra as convex hulls also.
The details are as follows.

We define a \emph{projective polytope} in $\oproj{d}$ to be the convex
hull of a finite set of points $V \subseteq \oproj{d}$,
all of which lie within the interior of some projective half-space.
If the set $V$ is minimal (that is, no point can be removed from $V$
without changing the convex hull), we call the points of $V$ the
\emph{vertices} of the projective polytope.

It is clear that traditional Euclidean polytopes correspond precisely to
projective polytopes whose vertices are all visible Euclidean points.
The situation regarding unbounded Euclidean polyhedra is more interesting.

Suppose that the set $V \subseteq \oproj{d}$ consists of a mixture of
visible Euclidean points and points at infinity, and (as usual) that
all of these points lie within the interior of some projective half-space.
Then the projective polytope with vertex set $V$ is in fact an
unbounded Euclidean polyhedron in visible space, ``closed off'' with
additional points at infinity.

\begin{figure}[htb]
    \centering
    \includegraphics[scale=0.9]{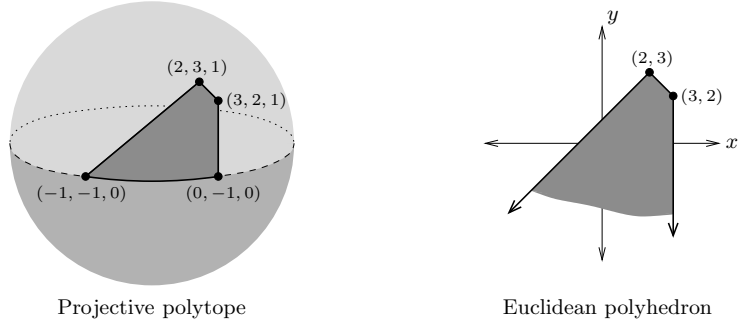}
    \caption{A 2-dimensional projective polytope with vertices at infinity}
    \label{fig-poly2}
\end{figure}

For example, Figure~\ref{fig-poly2} shows a projective
polytope in $\oproj{2}$, with two visible vertices
$(2,3,1)$ and $(3,2,1)$, and two infinite vertices
$(-1,-1,0)$ and $(0,-1,0)$.  The left-hand diagram shows this
projective polytope in $\oproj{2}$, and the right-hand diagram shows the
corresponding unbounded polyhedron in $\euc{2}$.

Conversely, let $P$ be an unbounded Euclidean polyhedron in $\euc{d}$
that lies within the interior of some Euclidean half-space.  Then
$P$ can be expressed as a projective polytope.  More precisely, there
is a projective polytope in $\oproj{d}$ that consists of the polyhedron $P$
in visible space ``closed off'' with an extra facet at infinity.
The vertices of this projective polytope are the vertices of
$P$ in visible space with some additional vertices at infinity.

\begin{figure}[htb]
    \centering
    \includegraphics[scale=0.9]{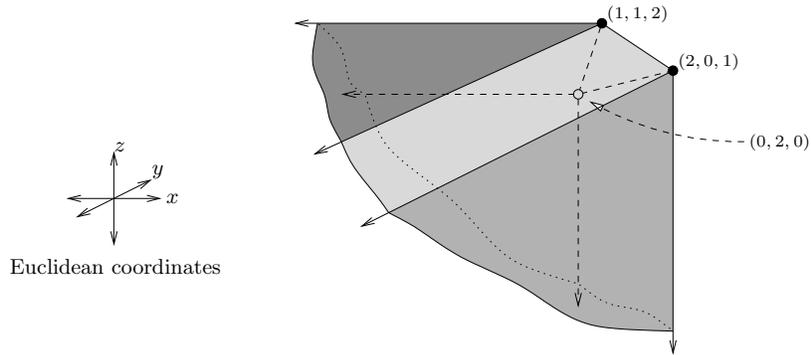}
    \caption{An unbounded 3-dimensional Euclidean polyhedron}
    \label{fig-poly3}
\end{figure}

For a 3-dimensional example, Figure~\ref{fig-poly3} shows a Euclidean
polyhedron in $\euc{3}$ with vertices $(2,0,1)$, $(1,1,2)$ and $(0,2,0)$,
and with facets that extend infinitely far in the negative $x$, $y$ and
$z$ directions.  The corresponding projective polytope has six vertices:
the visible vertices $(2,0,1,1)$, $(1,1,2,1)$ and $(0,2,0,1)$,
and the three infinite vertices $(-1,0,0,0)$, $(0,-1,0,0)$ and $(0,0,-1,0)$.

As in traditional geometry, any projective polytope can be expressed as
an intersection of finitely many projective half-spaces.  Moreover,
projective polytopes have the same combinatorial properties as Euclidean
polytopes, and they support the same combinatorial algorithms
(such as the double description method for vertex enumeration, which
we use later in this paper).\footnote{%
    This equivalence arises because any projective polytope can
    be transformed into a visible Euclidean polytope without changing its
    combinatorial properties.  In the spherical model of
    Figure~\ref{fig-proj-plane}, the transformation is essentially
    a rotation of the sphere.}

%% file: eff.tex
\section{Efficiency-equivalent polyhedra} \label{s-eff}

Recall that our focus is on generating all
\emph{efficient extreme outcomes} for a multiobjective linear program.
In other words, our aim is to generate all vertices of the outcome set
$Y^=$ (a polytope in $\euc{p}$) that are also efficient (that is, not
dominated by any other point in $Y^=$).

One difficulty is that the outcome set $Y^=$ can become
extremely complex even for relatively small problems
\cite{steuer05-regression,ziegler95}, and the efficient vertices
might only constitute a small fraction of the set of all vertices of $Y^=$.
Several authors address this by working with
\emph{efficiency-equivalent polyhedra}---polyhedra in $\euc{p}$
that are typically simpler in structure than $Y^=$, but whose efficient
sets precisely match the efficient outcome set $Y^=_E$.

There are two such polyhedra that feature prominently in the
literature: the bounded polytope $Y^\square$ introduced by Benson
\cite{benson98-outer}, and the unbounded polyhedron $Y^\leq$ first
described by Dauer and Saleh
\cite{dauer90-constructing,dauer92-representation},
both of which we describe shortly.\footnote{%
    The symbols $Y^=$ and $Y^\leq$ follow Benson's notation;
    Dauer and Saleh call these $Y$ and $\tilde{Y}$ instead.
    The symbol $Y^\square$ is our own notation;
    Benson simply calls this $Y$.
    We avoid the plain symbol $Y$ in this paper because it is used with
    several different meanings throughout the literature.}
The goal of this section is to understand the structure and
complexity of these two polyhedra.

We begin this section with a general result about efficiency-equivalent
polyhedra that justifies our focus on \emph{vertices} of these polyhedra.
We then describe the efficiency-equivalent polyhedra
$Y^\square$ and $Y^\leq$, and we reinterpret the unbounded polyhedron
$Y^\leq$ as a \emph{polytope} $T^\leq$ using oriented projective geometry.
The remainder of this section concentrates on the complexity of
the Euclidean polytope $Y^\square$ and the projective polytope $T^\leq$,
paying particular attention to the
number of ``unwanted'' \emph{non-efficient} vertices.
We find through our analysis that $T^\leq$ is far
simpler than $Y^\square$ (and therefore far superior)
in this respect.

\subsection{Definitions and general results}

Our first result shows that \emph{any} efficiency-equivalent polyhedron
for $Y^=$ can be used to generate the efficient extreme outcomes for our
multiobjective linear program.
Ours is a general result;
Benson \cite{benson98-outer} and Dauer and Saleh \cite{dauer90-constructing}
make similar observations for the specific polyhedra
$Y^\square$ and $Y^\leq$ respectively.

\begin{lemma} \label{l-eff-vertices}
    Let $Z \subseteq \euc{p}$ be any efficiency-equivalent polyhedron for
    $Y^=$.  Then the efficient vertices of $Z$ are precisely the
    efficient vertices of $Y^=$.  That is, the efficient vertices
    of $Z$ are precisely the efficient extreme outcomes for our
    multiobjective linear program.
\end{lemma}

\begin{proof}
    Let $\mathbf{v}$ be an efficient vertex of $Y^=$.
    By efficiency-equivalence we know that $\mathbf{v}$ lies in the
    efficient set $Z_E$.
    We aim to show that $\mathbf{v}$ is a vertex (and thus an
    efficient vertex) of $Z$.

    Let $F$ be the smallest-dimensional face of $Z$ containing $\mathbf{v}$.
    Because $Z_E$ is a union of faces of $Z$, every point in $F$ must be
    efficient.
    By efficiency-equivalence, the entire face $F$ must lie in $Y^=$ also.
    If $\mathbf{v}$ is not a vertex of $Z$ then
    $\mathbf{v}$ lies within the relative
    interior of the face $F$, and since $F \subseteq Y^=$ it follows that
    $\mathbf{v}$ cannot be a vertex of $Y^=$, a contradiction.

    Applying the same argument in reverse shows that every efficient
    vertex of $Z$ is an efficient vertex of $Y^=$, and the proof is
    complete.
\end{proof}

We proceed now to define the polyhedra $Y^\square$ and $Y^\leq$,
which will occupy our attention for the remainder of this section.
Both are ``extensions'' of the outcome set $Y^=$, in that they enlarge
the non-efficient portion of $Y^=$ to create a simpler structure.

\begin{defn} \label{d-leq-square}
    We define the Euclidean polyhedron $Y^\leq$ to be the set
    of all points in $\euc{p}$ that are dominated by any point in $Y^=$.
    That is,
    \[ Y^\leq = \left\{ \mathbf{z} \in \euc{p}\,\left|\,
        \mathbf{z} \leq \mathbf{y}\ \textrm{for some}
        \ \mathbf{y} \in Y^= \right.\right\}. \]
    Choose some point $\widehat{\mathbf{y}} \in \euc{p}$ that is
    strictly dominated by all of $Y^=$; in other words,
    $\widehat{\mathbf{y}} < \mathbf{y}$ for all $\mathbf{y} \in Y^=$.
    We then define the Euclidean polytope $Y^\square$ to be
    \[ Y^\square = \left\{ \mathbf{z} \in \euc{p}\,\left|\,
        \widehat{\mathbf{y}} \leq \mathbf{z} \leq \mathbf{y}
        \ \textrm{for some}\ \mathbf{y} \in Y^= \right.\right\}. \]
\end{defn}

\begin{figure}[htb]
    \centering
    \includegraphics[scale=0.9]{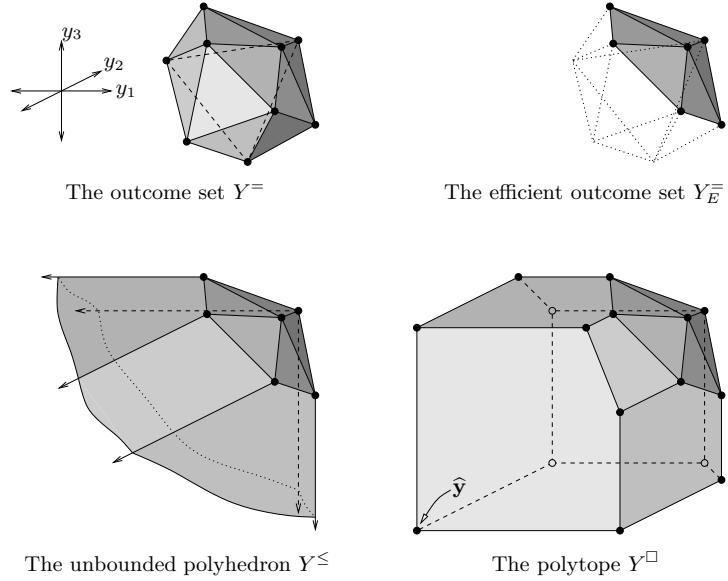}
    \caption{Constructing the efficiency-equivalent polyhedra
        $Y^\leq$ and $Y^\square$}
    \label{fig-eff-cube}
\end{figure}

Both of these structures are illustrated in Figure~\ref{fig-eff-cube}
(where we use $y_1,y_2,y_3$ for our three coordinate axes to show
that we are working in objective space).
It is clear that $Y^\leq$ is an unbounded polyhedron, and that it can
be expressed as the Minkowski sum
\begin{equation} \label{eqn-minkowski}
    Y^\leq = \left\{ \mathbf{y} + \mathbf{z}\,\left|\,
        \mathbf{y} \in Y^=\ \mbox{and}\ \mathbf{z} \leq \mathbf{0}
        \right.\right\}.
\end{equation}
Likewise, it is clear that $Y^\square$ is a polytope, since its points
are bounded below by $\widehat{\mathbf{y}}$ and above by $Y^=$.
Geometrically, $Y^\square$ truncates the unbounded polyhedron $Y^\leq$ by
intersecting it with a large cube that surrounds the original
outcome set $Y^=$.

The following crucial properties of $Y^\square$ and $Y^\leq$
are proven by Benson \cite{benson98-outer} and
Dauer and Saleh \cite{dauer90-constructing} respectively.

\begin{lemma} \label{l-known-equiv}
    Both $Y^\square$ and $Y^\leq$ are efficiency-equivalent polyhedra
    for $Y^=$.  In other words, the efficient sets
    $Y^\square_E$, $Y^\leq_E$ and $Y^=_E$ are identical.
\end{lemma}

\subsection{The complexity of $Y^\leq$ and $T^\leq$}

The unbounded polytope $Y^\leq$ is simpler than $Y^\square$ in structure,
with fewer vertices and a convenient Minkowski sum representation.
However, $Y^\leq$ is more difficult to work with algorithmically
because it is not a polytope---in particular, it cannot be
expressed as the convex hull of its vertices.
This makes it less amenable to traditional algorithms for
vertex enumeration and linear optimisation.

We now grant ourselves the best of both worlds by reinterpreting
the unbounded polyhedron $Y^\leq$ as a convex hull of vertices in
oriented projective $p$-space.

\begin{defn} \label{d-leq-proj}
    Let $T^\leq$ be the projective polytope in $\oproj{p}$
    formed as the convex hull of:
    \begin{enumerate}[(i)]
        \item \label{en-proj-t-visible}
        the vertices of the polyhedron $Y^\leq$ in
        visible Euclidean space;
        \item \label{en-proj-t-infinite}
        the points at infinity
        defined by the
        signed homogeneous coordinates
        $(-1,0,0,\ldots,0,0)$,
        $(0,-1,0,\ldots,0,0)$,
        $(0,0,-1,\ldots,0,0)$,
        \ldots,
        $(0,0,0,\ldots,-1,0)$
        (there are $p$ such points in total).
    \end{enumerate}
\end{defn}

In the following series of results, we show that the projective polytope
$T^\leq$ is well-defined, that it is truly a projective reinterpretation
of $Y^\leq$, and that Definition~\ref{d-leq-proj} gives the precise
vertex structure of $T^\leq$ in oriented projective $p$-space.

\begin{lemma} \label{l-leq-defined}
    The projective polytope $T^\leq$ is well-defined.
    That is, all of the points listed in (\ref{en-proj-t-visible}) and
    (\ref{en-proj-t-infinite}) above lie in the interior of some
    projective half-space.
\end{lemma}

\begin{proof}
    Because the Euclidean polyhedron $Y^\leq$ has finitely many vertices,
    there is some $\alpha \in \R$ that satisfies $\alpha > v_1+\ldots+v_p$
    for every Euclidean vertex $(v_1,\ldots,v_p)$ of $Y^\leq$.
    When translated into oriented projective $p$-space, it follows that
    every vertex of $Y^\leq$ lies in the interior of the projective
    half-space with signed homogeneous coordinates
    $(-1,\ldots,-1,\alpha)$.  The same is clearly true for each infinite
    point in list (\ref{en-proj-t-infinite}) above, and so the
    projective half-space $(-1,\ldots,-1,\alpha)$ is the half-space that
    we seek.
\end{proof}

\begin{lemma} \label{l-leq-visible}
    The visible portion of the projective polytope $T^\leq$
    is precisely the Euclidean polyhedron $Y^\leq$.
\end{lemma}

\begin{proof}
    We prove this by first showing that every visible Euclidean point of
    $T^\leq$ belongs to $Y^\leq$, and then that every point of
    the visible Euclidean polyhedron $Y^\leq$ belongs to $T^\leq$.
    \begin{itemize}
        \item Let $\mathbf{t} \in \oproj{p}$ be some visible Euclidean
        point in $T^\leq$.
        We first observe that
        any convex combination of points in list
        (\ref{en-proj-t-visible}) above is a visible Euclidean point in
        $Y^\leq$, and that any convex combination of points in list
        (\ref{en-proj-t-infinite}) above is a point at infinity with signed
        homogeneous coordinates $(-\beta_1,\ldots,-\beta_p,0)$ for some
        scalars $\beta_1,\ldots,\beta_p \geq 0$.
        Since $\mathbf{t}$ is a convex combination of points in
        \emph{both} lists, it follows that $\mathbf{t}$ has signed
        homogeneous coordinates of the form
        $(\mu y_1 - \lambda \beta_1, \ldots, \mu y_p - \lambda \beta_p, \mu)$,
        where $\mu,\lambda,\beta_1,\ldots,\beta_p \geq 0$
        and where $(y_1,\ldots,y_p)$ are the Euclidean coordinates of some
        point in $Y^\leq$.

        Because $\mathbf{t}$ is visible we must have $\mu > 0$, so we
        can rescale the coordinates of $\mathbf{t}$ to
        $(y_1 - \lambda \beta_1 / \mu, \ldots, y_p - \lambda \beta_p / \mu, 1)$.
        It is clear then that $\mathbf{t}$ is a visible point whose
        Euclidean coordinates are formed by adding the non-positive
        vector $-\lambda/\mu (\beta_1,\ldots,\beta_p)$ to some point in
        $Y^\leq$.  From equation~(\ref{eqn-minkowski})
        it follows that $\mathbf{t}$ is a visible Euclidean point
        in $Y^\leq$.

        \item Now let $\mathbf{y} \in \oproj{p}$ be a point of $Y^\leq$
        in visible Euclidean space with signed homogeneous coordinates
        $(y_1,\ldots,y_p,1)$.  Let $\mathbf{z}=(z_1,\ldots,z_p,1)$
        be an efficient point of
        $Y^\leq$ that dominates $\mathbf{y}$.  We can express
        the coordinates of $\mathbf{y}$ as
        $(z_1-\beta_1,\ldots,z_p-\beta_p,1)$
        for some scalars $\beta_1,\ldots,\beta_p \geq 0$,
        which shows that
        $\mathbf{y}$ is a convex combination of the efficient point
        $\mathbf{z}$ with the points at infinity in
        list~(\ref{en-proj-t-infinite}) above.

        Recalling that the efficient set $Y^\leq_E$ is a union
        of faces of $Y^\leq$, let $F$ be some
        efficient face of $Y^\leq$ containing $\mathbf{z}$.  Because
        $Y^\leq$ is efficiency-equivalent to the bounded polytope $Y^=$,
        the efficient face $F$ must likewise be bounded.
        Therefore $\mathbf{z}$ can be expressed as a convex combination
        of vertices of $F$, which are of course vertices of $Y^\leq$ also.
        This allows us to express $\mathbf{y}$ as a convex combination
        of points in lists~(\ref{en-proj-t-visible}) and
        (\ref{en-proj-t-infinite}), and so $\mathbf{y}$ belongs to the
        projective polytope $T^\leq$.
        \qedhere
    \end{itemize}
\end{proof}

\begin{lemma} \label{l-leq-vertices}
    The vertices of $T^\leq$ are precisely the points listed
    in Definition~\ref{d-leq-proj}.
    In other words, no point can be removed from these lists
    without changing the convex hull.
\end{lemma}

\begin{proof}
    To prove this we must show that no point in lists
    (\ref{en-proj-t-visible}) and (\ref{en-proj-t-infinite}) from
    Definition~\ref{d-leq-proj} can be
    expressed as a convex combination of the others.
    Suppose that list~(\ref{en-proj-t-visible}) contains the points
    $\mathbf{v}_1,\ldots,\mathbf{v}_t$ (these are the vertices of
    $Y^\leq$), and that list~(\ref{en-proj-t-infinite}) contains the
    points $\mathbf{f}_1,\ldots,\mathbf{f}_p$ (these are all points at
    infinity).
    \begin{itemize}
        \item No $\mathbf{f}_i$ can be expressed as a convex combination
        of the other points in these lists, since each other point
        introduces some non-zero coordinate that is zero in $\mathbf{f}_i$.

        \item Suppose that some $\mathbf{v}_i$ can be expressed as
        the convex combination $\sum_{j \neq i} \alpha_j \mathbf{v}_j +
        \sum_{k} \gamma_k \mathbf{f}_k$ where each $\alpha_j,\gamma_k \geq 0$
        (as usual we work in homogeneous coordinates).
        Some $\alpha_j$ must be non-zero, since otherwise $\mathbf{v}_i$
        would be a point at infinity.
        Some $\gamma_k$ must also be non-zero, since we cannot express a
        vertex of $Y^\leq$ as a convex combination of other vertices.

        It follows that we can express $\mathbf{v}_i$ in the form
        $(y_1-\beta_1,\ldots,y_p-\beta_p,1)$, where
        $\mathbf{y}=(y_1,\ldots,y_p,1)$ is a convex combination of
        vertices of $Y^\leq$ (in particular, $\mathbf{y} \in Y^\leq$),
        and where $\beta_1,\ldots,\beta_p \geq 0$ with at least
        one $\beta_k$ strictly positive.
        This places $\mathbf{v}_i$ in the interior of the line
        segment joining $\mathbf{y}$ with the point
        $(y_1-2\beta_1,\ldots,y_p-2\beta_p,1)$.  However, both endpoints of
        this line segment lie in the polyhedron $Y^\leq$, which is
        impossible since $\mathbf{v}_i$ is a vertex of $Y^\leq$.
        \qedhere
    \end{itemize}
\end{proof}

We can draw some important conclusions from these results.
The first is that $T^\leq$ is precisely the
efficiency-equivalent polyhedron $Y^\leq$, recast as a \emph{polytope} using
oriented projective geometry.  Furthermore, Lemma~\ref{l-leq-vertices}
gives an insight into the \emph{complexity} of this polytope, as seen in
the following result.

\begin{corollary} \label{c-leq-complexity}
    The visible Euclidean vertices of $T^\leq$ are precisely the
    efficient extreme outcomes for our multiobjective linear program.
    That is, the only vertices of $T^\leq$ that are not efficient
    extreme outcomes are the $p$ vertices at infinity.
\end{corollary}

\begin{proof}
    From Lemma~\ref{l-leq-vertices}, the visible Euclidean vertices of
    $T^\leq$ are precisely the vertices of the unbounded polyhedron
    $Y^\leq \subseteq \euc{p}$.
    Moreover, Lemma~\ref{l-eff-vertices} shows that the
    \emph{efficient} vertices of $Y^\leq$ are precisely the
    efficient extreme outcomes.

    Suppose that $Y^\leq$ has some \emph{non-efficient}
    vertex $\mathbf{v} \in Y^\leq$.
    Then $\mathbf{v}$ must be dominated by some other
    point $\mathbf{x} \in Y^\leq$, so that the difference
    $\mathbf{d} = \mathbf{x} - \mathbf{v}$ is non-negative and non-zero.
    Now construct the point $\mathbf{x}' = \mathbf{v} - \mathbf{d}$.
    It is clear that $\mathbf{x}' \in Y^\leq$ and that
    $\mathbf{v}$ lies strictly inside the line segment
    $[\mathbf{x},\mathbf{x}'] \subseteq Y^\leq$, contradicting the
    claim that $\mathbf{v}$ is a vertex of $Y^\leq$.
\end{proof}

It follows from Corollary~\ref{c-leq-complexity} that
the projective polytope $T^\leq$ gives an extremely compact
representation of the efficient extreme outcomes---the number of
additional ``unwanted'' vertices grows linearly with the number of
objectives $p$, and does not depend at all on the problem size
$(n,m)$.

\subsection{The complexity of $Y^\square$}

We turn now to study the complexity of the Euclidean polytope
$Y^\square$.  We begin with a complete categorisation of vertices of
$Y^\square$, and we follow with a discussion of what this means for
best-case and worst-case complexity.
Throughout the remainder of this section we work entirely in Euclidean
coordinates.

\begin{defn} \label{d-vset}
    For each set of indices $S \subseteq \{1,\ldots,p\}$, we define the
    set of points $V(S) \subseteq \euc{p}$ as follows.  Denote the
    $p$ objectives by $y_1,\ldots,y_p$, and consider the smaller
    multiobjective linear program obtained by considering only the
    objectives $\{y_i\,|\,i \in S\}$ (in other words, by deleting
    the $i$th row of the objective matrix for each $i \notin S$).
    This smaller problem yields a smaller set
    of efficient extreme outcomes in the space $\euc{|S|}$.

    We extend these efficient extreme outcomes to points in $\euc{p}$
    as follows.  For each deleted objective $i \notin S$, we fix the
    coordinate $y_i$ to the value $\widehat{y}_i$
    (using the point $\widehat{\mathbf{y}} \in \euc{p}$ that we describe
    in Definition~\ref{d-leq-square}).
    The resulting set of points in $\euc{p}$ is denoted $V(S)$.

    If $S$ is the empty set then there is no linear program to solve,
    and we simply define $V(\emptyset) = \{\widehat{\mathbf{y}}\}$.
\end{defn}

\begin{figure}[htb]
    \centering
    \includegraphics[scale=0.9]{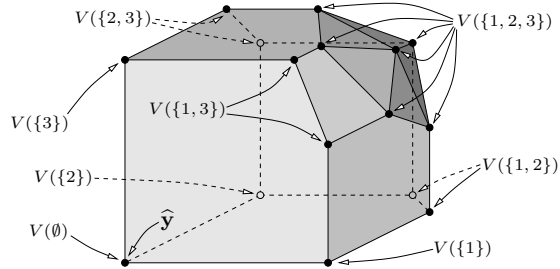}
    \caption{The sets $V(S)$ in relation to the polytope $Y^\square$}
    \label{fig-eff-v}
\end{figure}

It is worth noting that for the full set $S = \{1,\ldots,p\}$,
the points in $V(S)$ are simply the efficient extreme outcomes for our
original (full-dimensional) multiobjective linear program.

The various sets $V(S)$ are illustrated in Figure~\ref{fig-eff-v},
which continues the earlier example of Figure~\ref{fig-eff-cube}.
We see in this diagram that the sets $V(S)$ together describe all
of the vertices of $Y^\square$.  This is no accident,
as shown by the following result.

\begin{theorem} \label{t-vert-square}
    The sets $V(S)$ are pairwise disjoint; that is,
    if $S \neq S'$ then $V(S) \cap V(S') = \emptyset$.
    Moreover, the union of all such sets $V(S)$ over all
    $S \subseteq \{1,\ldots,p\}$ gives the complete vertex set for the
    polytope $Y^\square$.
\end{theorem}

\begin{proof}
    Pairwise disjointness is simple to prove.
    Because $\widehat{\mathbf{y}}$ is strictly dominated by all
    $\mathbf{y} \in Y^=$, we find that for any
    $\mathbf{v} = (v_1,\ldots,v_p) \in V(S)$
    we have
    $v_i = \widehat{y}_i$ if and only if $i \in S$.
    It is then clear that $\mathbf{v}$ cannot belong
    to both $V(S)$ and $V(S')$ for distinct sets $S \neq S'$.

    It remains to show that every element of each set $V(S)$ is a
    vertex of $Y^\square$, and that every vertex of $Y^\square$ can be
    expressed as an element of some $V(S)$.
    For this it is convenient to define the linear projection
    $\pi_S \co \euc{p} \to \euc{|S|}$ that
    simply deletes all coordinates $y_i$ for $i \notin S$.
    For example, if $S = \{1,\ldots,p\}$ then $\pi_S$ is the identity map,
    and if $S = \{i\}$ then $\pi_S$ merely extracts the $i$th coordinate.

    From Definition~\ref{d-vset} we have $\mathbf{v} \in V(S)$
    if and only if (i)~$v_i = \widehat{y}_i$ for each $i \notin S$,
    and (ii)~the projection $\pi_S(\mathbf{v})$ is an
    efficient vertex of the polytope $\pi_S(Y^=)$.
    However, applying Lemma~\ref{l-known-equiv} to the smaller
    multiobjective linear program with objectives $\{y_i\,|\,i \in S\}$
    shows that $\pi_S(Y^\square)$ is efficiency-equivalent to $\pi_S(Y^=)$.
    This means that $\mathbf{v} \in V(S)$ if and only if
    (i)~$v_i = \widehat{y}_i$ for each $i \notin S$, and
    ($\mathrm{ii}'$)~$\pi_S(\mathbf{v})$
    is an efficient vertex of $\pi_S(Y^\square)$.

    We can now conclude the proof as follows:
    \begin{itemize}
        \item
        Suppose that $\mathbf{v} \in V(S)$ for some set $S \neq \emptyset$.
        We aim to show that $\mathbf{v}$ is a vertex of $Y^\square$.
        If not, we can express $\mathbf{v} = \frac12(\mathbf{u}+\mathbf{w})$
        for some distinct $\mathbf{u},\mathbf{w} \in Y^\square$.
        We proceed to derive a contradiction by showing that
        $\mathbf{v} = \mathbf{u} = \mathbf{w}$.
        \begin{itemize}
            \item
            For any $i \notin S$ we have $v_i = \widehat{y}_i$ by
            definition of $V(S)$, and $u_i,w_i \geq \widehat{y}_i$ by
            definition of $Y^\square$.  Since
            $\mathbf{v} = \frac12(\mathbf{u}+\mathbf{w})$ it follows
            that $v_i = u_i = w_i$.

            \item
            By linearity, $\pi_S(\mathbf{v}) =
            \frac12\left(\pi_S(\mathbf{u})+\pi_S(\mathbf{w})\right)$
            where $\pi_S(\mathbf{u}),\pi_S(\mathbf{w}) \in \pi_S(Y^\square)$.
            However, because $\mathbf{v} \in V(S)$ we know that
            $\pi_S(\mathbf{v})$ is a vertex of $\pi_S(Y^\square)$,
            which is only possible if
            $\pi_S(\mathbf{v}) = \pi_S(\mathbf{u}) = \pi_S(\mathbf{w})$.
            That is, $v_i = u_i = w_i$ for all $i \in S$.
        \end{itemize}
        This establishes $\mathbf{v} = \mathbf{u} = \mathbf{w}$,
        and we conclude that $\mathbf{v}$ must indeed be a vertex
        of $Y^\square$.

        We handle the case $S = \emptyset$ by
        observing that $V(\emptyset)$ contains the single point
        $\widehat{\mathbf{y}}$, and that $\widehat{\mathbf{y}}$ is a vertex
        of $Y^\square$ because it uniquely minimises the linear
        functional $\mathbf{y} \mapsto y_1+\ldots+y_p$.

        \item
        Now let $\mathbf{v}$ be some vertex of $Y^\square$, and let
        $S = \{i\,|\,v_i \neq \widehat{y}_i\}$.
        We aim to show that $\mathbf{v} \in V(S)$.
        If $S = \emptyset$ then $\mathbf{v} = \widehat{\mathbf{y}}
        \in V(\emptyset)$, so assume that $S$ is non-empty.
        We already have $v_i = \widehat{y}_i$ for each $i \notin S$, so
        we only need to show that $\pi_S(\mathbf{v})$ is an efficient
        vertex of $\pi_S(Y^\square)$.

        \begin{itemize}
            \item If $\pi_S(\mathbf{v})$ is not a vertex of
            $\pi_S(Y^\square)$, there must be some
            $\mathbf{u},\mathbf{w} \in Y^\square$ for which
            $\pi_S(\mathbf{v}) =
            \frac12\left(\pi_S(\mathbf{u})+\pi_S(\mathbf{w})\right)$
            but $\pi_S(\mathbf{u}) \neq \pi_S(\mathbf{w})$.
            Construct the points $\mathbf{u}',\mathbf{w}'$ from
            $\mathbf{u},\mathbf{w}$ respectively by replacing each $i$th
            coordinate with $\widehat{y}_i$ for all $i \notin S$.

            It is clear that $\mathbf{u}',\mathbf{w}' \in Y^\square$
            (since $\widehat{\mathbf{y}} \leq \mathbf{u}' \leq \mathbf{u}$
            and $\widehat{\mathbf{y}} \leq \mathbf{w}' \leq \mathbf{w}$).
            Moreover, for $i \notin S$ we have $v_i = u'_i = w'_i$,
            and for $i \in S$ we have $v_i = \frac12(u'_i + w'_i)$.
            Finally, we note that $\pi_S(\mathbf{u}) \neq \pi_S(\mathbf{w})$
            implies that $\mathbf{u}' \neq \mathbf{w}'$.
            We therefore have $\mathbf{v} = \frac12(\mathbf{u'}+\mathbf{w'})$
            but $\mathbf{u'} \neq \mathbf{w'}$, contradicting the
            assumption that $\mathbf{v}$ is a vertex of $Y^\square$.

            \item If $\pi_S(\mathbf{v})$ is not efficient in
            $\pi_S(Y^\square)$, there must be some
            $\mathbf{u} \in Y^\square$ for which
            $\pi_S(\mathbf{u}) \geq \pi_S(\mathbf{v})$ but
            $\pi_S(\mathbf{u}) \neq \pi_S(\mathbf{v})$.
            Again, construct $\mathbf{u}'$ from $\mathbf{u}$ by replacing
            the $i$th coordinate with $\widehat{y}_i$ for each $i \notin S$.
            This time we have $\mathbf{u}' \in Y^\square$,
            $\mathbf{u}' \neq \mathbf{v}$, and the
            relations $\widehat{y}_i = v_i = u'_i$ for $i \notin S$ and
            $\widehat{y}_i < v_i \leq u'_i$ for $i \in S$.

            Let $\mathbf{z} =
            \mathbf{v} - \varepsilon (\mathbf{u}'-\mathbf{v})$ for
            some small $\varepsilon > 0$ that ensures
            $\widehat{\mathbf{y}} \leq \mathbf{z}$ (by the relations
            above, such an $\varepsilon$ exists).
            Then $\mathbf{z} \in Y^\square$ and
            $\mathbf{v}$ is in the interior of the line segment
            $[\mathbf{z},\mathbf{u}']$, again contradicting
            the assumption that $\mathbf{v}$ is a vertex of $Y^\square$.
        \end{itemize}
        Together these observations show that $\pi_S(\mathbf{v})$ is an
        efficient vertex of $\pi_S(Y^\square)$, whereupon the proof is
        complete.
        \qedhere
    \end{itemize}
\end{proof}

Theorem~\ref{t-vert-square} has important implications for the complexity of
the polytope $Y^\square$.  Recall that all of the efficient vertices of
$Y^\square$ are held in the single set $V(\{1,\ldots,p\})$.
This means that for every \emph{proper} subset $S \subset \{1,\ldots,p\}$,
the corresponding $V(S)$ is filled with \emph{non-efficient} vertices of
$Y^\square$.  This can yield an extremely large number of non-efficient
vertices---each $V(S)$ describes the efficient extreme outcomes for some
other multiobjective linear program with up to $p-1$ objectives, and
we already know that the number of such outcomes can grow very large
even for relatively small values of $p$, $n$ and $m$
\cite{steuer05-regression,ziegler95}.

This argument shows intuitively that $Y^\square$ can have an extremely large
number of non-efficient vertices.  We follow now with two corollaries of
Theorem~\ref{t-vert-square} that back this up with explicit formulae.
The first result shows that, even in the \emph{best case} scenario where
$Y^\square$ is as simple as possible, we still obtain an exponentially
large number of non-efficient vertices.

\begin{corollary} \label{c-square-best}
    The polytope $Y^\square$ always has at least $2^p-1$ non-efficient
    vertices; that is, $2^p-1$ vertices that are not efficient extreme
    outcomes for our original multiobjective linear program.
\end{corollary}

\begin{proof}
    As noted above, every proper subset
    $S \subset \{1,\ldots,p\}$ yields a corresponding $V(S)$
    that consists entirely of non-efficient vertices.
    There are $2^p-1$ proper subsets of $\{1,\ldots,p\}$,
    and so we have at least $2^p-1$ non-efficient vertices in total.
\end{proof}

The corollary above shows that the growth rate of non-efficient
vertices relative to the number of objectives $p$ is unavoidably exponential.
In our second corollary, we show that the number of non-efficient
vertices can also grow at a severe rate relative to the
problem size $(n,m)$.
This is an asymptotic result, and we use the
standard notation for complexity whereby
$O(\cdot)$ and $\Omega(\cdot)$ denote asymptotic upper bounds
and lower bounds respectively.

\begin{corollary} \label{c-square-worst}
    Suppose we fix $p$ and allow the problem size $(n,m)$ to vary.
    Then there are multiobjective linear programs for which
    the number of non-efficient vertices of
    $Y^\square$ grows at least as fast as
    $\Omega(n^{\lfloor (p-1)/2 \rfloor})$.
\end{corollary}

\begin{proof}
    We establish this result using \emph{dual cyclic polytopes}.
    For any integers $2 \leq d < k$, the $(d,k)$ dual cyclic polytope
    is a $d$-dimensional polytope with precisely $k$ facets and
    \begin{equation} \label{eqn-cyclic}
        \binom{k - \lfloor \frac{d+1}{2} \rfloor}{k - d} +
        \binom{k - \lfloor \frac{d+2}{2} \rfloor}{k - d}
    \end{equation}
    vertices.  For further information on dual cyclic polytopes the reader is
    referred to a standard reference such as Gr\"unbaum \cite{grunbaum03}.

    For us, the important fact about dual cyclic polytopes is that
    for a fixed dimension $d$, the number of vertices grows at least as
    fast as $\Omega(k^{\lfloor d/2 \rfloor})$.  This is because the sum
    of binomial coefficients in equation~(\ref{eqn-cyclic}) can be
    rewritten as a $\lfloor d/2 \rfloor$-degree polynomial in $k$
    (taking separate cases for odd and even $d$).

    Our strategy then is as follows.  Suppose we are given some number
    of objectives $p$ and some number of variables $n$.
    We assume here that $3 \leq p \leq n$ (since we are interested in the
    growth rate as $n$ becomes large, and for $p \leq 2$ the result
    is trivial).  Let $H$ be the hyperplane in $\euc{p}$ defined by the
    equation $y_1 + \ldots + y_p = 1$.  We now embed a $(p-1,n)$ dual
    cyclic polytope $P$ within the hyperplane $H$ (a straightforward
    procedure since $\dim H = p-1$).

    The polytope $P \subseteq \euc{p}$ can be expressed as the outcome set
    $Y^=$ for
    some multiobjective linear program with $p$ objectives and $n$ variables
    (its $n$ facets arise from the inequalities $\mathbf{x} \geq 0$ in the
    solution space $\euc{n}$).
    Although $P$ is not full-dimensional in $\euc{p}$, this is just
    a convenience to keep our arithmetic simple.
    If a full-dimensional example is ever required then we simply extend
    $P$ to a $p$-dimensional cone; the analysis becomes messier but the
    final growth rate of $\Omega(n^{\lfloor (p-1)/2 \rfloor})$ remains
    the same.

    We can now make the following observation.
    \begin{quote}
        \textbf{Claim:}
        \emph{Every vertex of $Y^=$ (that is, $P$) gives an efficient extreme
        outcome for some smaller multiobjective linear program that uses
        only $p-1$ of our $p$ objectives.}
    \end{quote}
    We will prove this claim shortly.  In the meantime, we examine its
    consequences.  By Theorem~\ref{t-vert-square}, it follows that each
    vertex of $P$ corresponds to some element of $V(S)$ where $S$ is some
    $(p-1)$-element subset of $\{1,\ldots,p\}$.  Moreover, because $P$
    lies entirely within the hyperplane $\sum y_i = 1$, any two vertices of
    $P$ must differ in \emph{at least two} coordinates.  It follows that
    any two distinct vertices of $P$ must correspond to either two distinct
    elements of the same set $V(S)$, or else elements of two distinct sets
    $V(S)$, $V(S')$.

    As a result, the sum of $|V(S)|$ over all $(p-1)$-element sets $S$
    is at least the number of vertices of $P$.  However, each such $S$ is a
    proper subset of $\{1,\ldots,p\}$, and we can therefore deduce that
    the total number of
    non-efficient vertices of $Y^\square$ is at least the number of
    vertices of $P$.  The number of vertices of $P$ is described by
    equation~(\ref{eqn-cyclic}), which we noted earlier grows at least as
    fast as
    $\Omega(k^{\lfloor d/2 \rfloor}) = \Omega(n^{\lfloor (p-1)/2 \rfloor})$.
    Thus our final result is established.

    All that remains then is to prove the claim above.
    Let $\mathbf{v}$ be any vertex of $Y^=$.  It follows that
    $\mathbf{v}$ is the unique intersection of $Y^=$ with some supporting
    hyperplane $H \subseteq \euc{p}$.  Suppose that $H$ is defined by
    the equation $\alpha_1 y_1 + \ldots + \alpha_p y_p = 0$, where $Y^=$ lies
    in the closed half-space $\alpha_1 y_1 + \ldots + \alpha_p y_p \leq 0$.

    We can ``tilt'' this hyperplane as follows.
    Let $\alpha_\mathrm{min} = \min\{\alpha_i\}$, and consider instead the
    hyperplane $H' \subseteq \euc{p}$ defined by the equation
    $(\alpha_1-\alpha_\mathrm{min}) y_1 + \ldots +
    (\alpha_p-\alpha_\mathrm{min}) y_p = 0$.
    Because every point in $Y^=$ satisfies $\sum y_i = 1$, we find that
    $H'$ is also a supporting hyperplane for $Y^=$,
    that $Y^=$ lies in the closed half-space
    $(\alpha_1-\alpha_\mathrm{min}) y_1 + \ldots +
    (\alpha_p-\alpha_\mathrm{min}) y_p \leq 0$,
    and that $\mathbf{v}$ is the unique intersection of $Y^=$ with $H'$.

    More importantly, at least one of the coefficients
    $(\alpha_i-\alpha_\mathrm{min})$ is zero, and all of the remaining
    coefficients of $H'$ are non-negative.  It follows that
    $(\alpha_1-\alpha_\mathrm{min}) y_1 + \ldots +
    (\alpha_p-\alpha_\mathrm{min}) y_p$ is a linear functional with
    non-negative coefficients in at most $p-1$ objectives, and that
    $\mathbf{v}$ is the unique maximum for this linear functional
    in the outcome set $Y^=$.  As a result, $\mathbf{v}$ indeed gives
    an efficient extreme outcome for some smaller multiobjective
    linear program using only $p-1$ of our $p$ objectives.
\end{proof}

%% file: outer.tex
\section{The outer approximation algorithm} \label{s-outer}

In this section we illustrate the power of our techniques
by applying them to the outer approximation
algorithm, which enumerates efficient extreme outcomes
for a multiobjective linear program.
Introduced by Benson \cite{benson98-outer}, this algorithm iteratively
constructs the efficiency-equivalent polytope $Y^\square$ using
methods from linear programming, univariate search techniques and
polytope vertex enumeration algorithms.
The algorithm is designed to work with polytopes, not unbounded
polyhedra, which explains the choice of $Y^\square$ as the final target.

We begin in Section~\ref{s-outer-benson} by describing the original outer
approximation algorithm in detail.  In Section~\ref{s-outer-new} we
develop a new outer approximation algorithm by
working in oriented projective geometry,
which allows us to replace the target polytope $Y^\square$ with the
much simpler projective polytope $T^\leq$.
In Section~\ref{s-outer-complexity} we examine how
these changes affect the runtime complexity.

\subsection{The original algorithm} \label{s-outer-benson}

\begin{figure}[htb]
    \centering
    \includegraphics[scale=0.3]{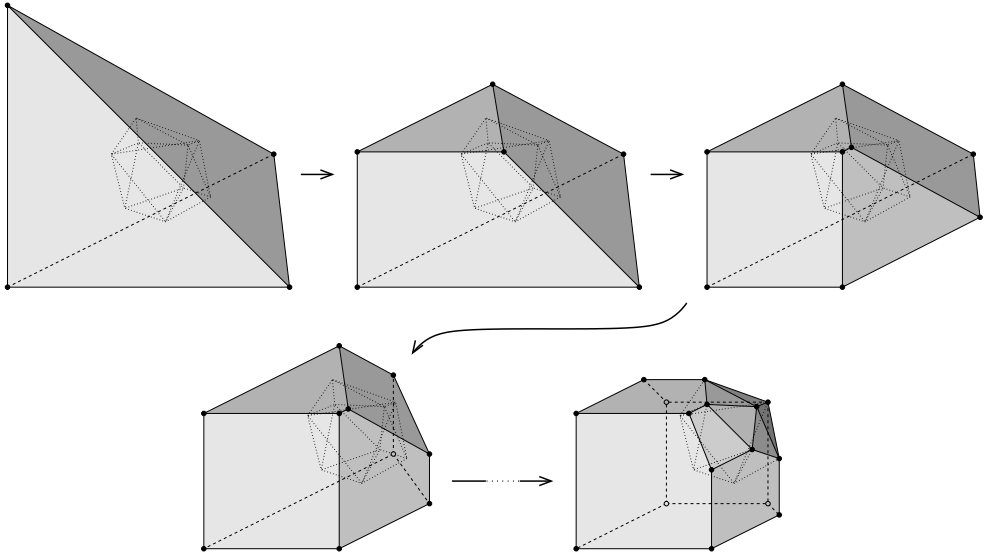}
    \caption{Stepping through the original outer approximation algorithm}
    \label{fig-outer}
\end{figure}

The algorithm originally described by Benson begins with a
simplex in $\euc{p}$ that completely surrounds the outcome set $Y^=$,
and then successively truncates this simplex along hyperplanes until it
is eventually whittled down to the target polytope $Y^\square$.
This overall procedure is illustrated in Figure~\ref{fig-outer}
(which continues the earlier example from Figure~\ref{fig-eff-cube}).
The details are as follows.

\begin{algorithm}[Outer approximation algorithm] \label{a-outer}
    This algorithm builds $Y^\square$ by constructing a series of
    polytopes $P_0 \supseteq P_1 \supseteq \ldots \supseteq P_k = Y^\square$
    in $\euc{p}$ (where the number of polytopes $k$ is not known in advance).
    Each polytope is stored
    using two representations: as a convex hull of vertices in
    $\euc{p}$, and as an intersection of half-spaces in $\euc{p}$.

    To initialise the algorithm,
    choose an arbitrary point $\overline{\mathbf{y}}$
    in the interior of $Y^\square$,
    and set the initial polytope $P_0$ to be a simplex with vertex set
    $\{\widehat{\mathbf{y}}$,
    $\widehat{\mathbf{y}} + \alpha \mathbf{e}_1$,
    \ldots,
    $\widehat{\mathbf{y}} + \alpha \mathbf{e}_p\}$.
    Here $\widehat{\mathbf{y}}$ is the point chosen in
    Definition~\ref{d-leq-square},
    each $\mathbf{e}_i$ represents the $i$th unit vector in $\euc{p}$,
    and the constant $\alpha$ is chosen so that the simplex contains the
    entire outcome set $Y^=$.

    Now iteratively construct each polytope $P_{i+1}$ from the previous
    polytope $P_i$ as follows:
    \begin{enumerate}[Step 1.]
        \item \label{en-outer-search}
        Search through the vertices of $P_i$ for one that
        does not belong to the target polytope $Y^\square$.
        If no such vertex exists (which indicates that $P_i = Y^\square$),
        then output the efficient vertices of $P_i$ and
        terminate the algorithm.

        \item \label{en-outer-bdry}
        Otherwise let $\mathbf{v}$ be such a vertex,
        and compute the unique point $\mathbf{x} \in \euc{p}$ at which
        the line segment joining $\mathbf{v}$ with $\overline{\mathbf{y}}$
        crosses the boundary of the target polytope $Y^\square$.

        \item \label{en-outer-half}
        Compute a half-space $H \subseteq \euc{p}$ that contains
        the target polytope $Y^\square$, and that has the point
        $\mathbf{x}$ in its bounding hyperplane.

        \item \label{en-outer-hull}
        The next polytope is now $P_{i+1} = P_i \cap H$.
        Compute a half-space representation of $P_{i+1}$ (which is
        trivial) and a convex hull representation (which is not).
    \end{enumerate}
\end{algorithm}

There are significant details to be filled in for each of these steps,
for which we refer the reader to Benson's
original paper \cite{benson98-outer}.  In brief:
\begin{itemize}
    \item In the initialisation phase, both the constant $\alpha$ and the
    interior point $\overline{\mathbf{y}} \in Y^\square$ can be found
    by solving a linear program,
    and the half-space representation of the initial simplex $P_0$ is
    then simple to compute.

    \item In Step~\ref{en-outer-search} of the iteration,
    we can test for membership in $Y^\square$ by testing the feasibility of
    another linear program.  The boundary point $\mathbf{x}$ in
    Step~\ref{en-outer-bdry} is found by combining feasibility
    tests with standard univariate search techniques (Benson uses
    a simple bisection search), and the half-space $H$ in
    Step~\ref{en-outer-half} uses the solution to a dual linear program.

    \item In the final output phase, the efficient vertices of
    $P_i=Y^\square$ are precisely those vertices that are
    \emph{strictly} dominated by $\widehat{\mathbf{y}}$.
\end{itemize}

The only step not covered in the list above is Step~\ref{en-outer-hull},
which requires us to
compute half-space and convex hull representations for the new
polytope $P_{i+1}$.  The half-space representation is simple: just append
$H$ to the list of half-spaces for the previous polytope $P_i$.  The
convex hull representation is more complex, and requires iterative vertex
enumeration techniques.  These are typically based on the
\emph{double description method} of Motzkin et~al.\ \cite{motzkin53-dd},
though the central idea has been rediscovered by many different authors
and goes under a variety of names.  See
Chen et~al.\ \cite{chen91-online} and
Horst et~al.\ \cite{horst88-vertices} for alternative treatments.

The fundamental principle behind the double description method is the
following result:

\begin{theorem} \label{t-dd}
    Let $P \subseteq \euc{p}$ be a polytope, and let
    $H \subseteq \euc{p}$ be a half-space.  Then the vertices of the
    polytope $P \cap H$ can be obtained from the vertices of $P$ using
    the following procedure.

    First partition the vertices of $P$ into three sets $S_+$,
    $S_0$ and $S_-$, according to whether each vertex lies
    within the interior of $H$, on the bounding hyperplane of $H$,
    or external to $H$ respectively.

    Each vertex in $S_0$ or $S_+$ then becomes a vertex of $P \cap H$.
    Furthermore, for each pair of vertices
    $\mathbf{u} \in S_+$ and $\mathbf{v} \in S_-$ that are adjacent in $P$,
    the point where the line segment $[\mathbf{u},\mathbf{v}]$
    cuts the bounding hyperplane of $H$ also becomes a vertex of $P \cap H$.
    The polytope $P \cap H$ has no other vertices besides those
    described here.
\end{theorem}

\begin{figure}[htb]
    \centering
    \includegraphics[scale=0.9]{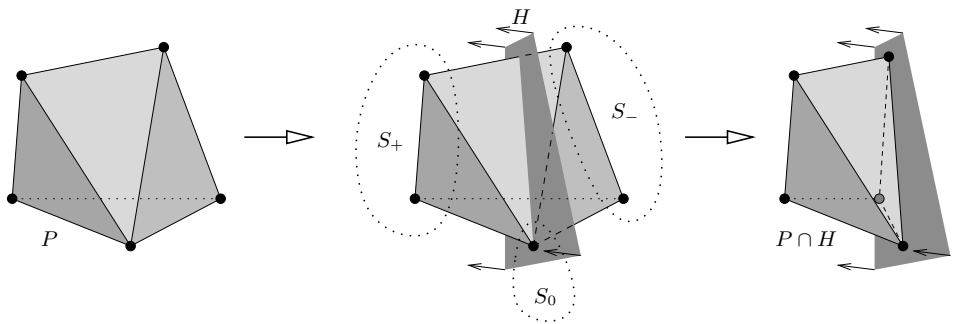}
    \caption{Intersecting the polytope $P$ with the half-space $H$}
    \label{fig-dd}
\end{figure}

This procedure is illustrated in Figure~\ref{fig-dd}.  The only
complication of Theorem~\ref{t-dd} is determining which pairs of
vertices $\mathbf{u} \in S_+$ and $\mathbf{v} \in S_-$ are adjacent
(i.e., joined by an edge) in the original polytope $P$.
There are many different adjacency tests in the literature;
in this paper we use the \emph{combinatorial adjacency test},
which is found to be fast and simple in practice \cite{fukuda96-doubledesc},
and which translates well to the setting of oriented projective geometry
(as seen later in Theorem~\ref{t-dd-proj}).

\begin{lemma} \label{l-adj}
    Let $P \subseteq \euc{p}$ be a polytope with half-space
    representation $P = H_1 \cap H_2 \cap \ldots \cap H_k$,
    where each $H_k$ is a half-space in $\euc{p}$.
    Then two vertices $\mathbf{u},\mathbf{v}$ of $P$ are adjacent if and
    only if there is no other vertex $\mathbf{z}$ of $P$ with the
    following property:
    \begin{itemize}
        \item
        For each half-space $H_i$, if both $\mathbf{u}$ and $\mathbf{v}$
        lie on the bounding hyperplane of $H_i$ then $\mathbf{z}$ lies
        on the bounding hyperplane also.
    \end{itemize}
\end{lemma}

For proofs of Theorem~\ref{t-dd} and Lemma~\ref{l-adj} as well as
further details on the double description
method, the reader is referred to Fukuda and Prodon \cite{fukuda96-doubledesc}.

\subsection{A new algorithm} \label{s-outer-new}

We now apply the results of Section~\ref{s-eff} to the outer
approximation algorithm, yielding a new algorithm with significant
benefits over the original.
Our strategy is to work in the oriented
projective space $\oproj{p}$ instead of the Euclidean space $\euc{p}$.
This allows us to iteratively construct the projective polytope $T^\leq$
instead of the Euclidean polytope $Y^\square$, thereby removing an
exponential number of non-efficient vertices that would otherwise
clutter up our computations.  We discuss the running time benefits
in more detail in Section~\ref{s-outer-complexity}.

As with the original algorithm, the final output for our new algorithm is the
set of all efficient extreme outcomes for our multiobjective linear program.
Because we work in oriented projective space, we perform all
calculations using signed homogeneous coordinates.
The details are as follows.

\begin{algorithm}[New outer approximation algorithm] \label{a-imp}
    In this algorithm we build the projective polytope
    $T^\leq$ by constructing a series of projective
    polytopes $Q_0 \supseteq Q_1 \supseteq \ldots \supseteq Q_k = T^\leq$
    in $\oproj{p}$.
    As usual each polytope is stored
    using two representations: as a convex hull of vertices in
    $\oproj{p}$, and as an intersection of projective half-spaces in
    $\oproj{p}$.

    To initialise the algorithm, construct the point
    $\mathbf{y}^\mathrm{max}$ with signed homogeneous coordinates
    $(y^\mathrm{max}_1,\ldots,y^\mathrm{max}_p,1)$,
    where each $y^\mathrm{max}_i$ individually maximises the
    $i$th objective.
    Choose an arbitrary point $\overline{\mathbf{y}}$
    in the interior of $T^\leq$, and set the
    initial polytope $Q_0$ to the simplex with
    one visible Euclidean vertex $\mathbf{y}^\mathrm{max}$
    and $p$ vertices at infinity with signed homogeneous coordinates
    $(-1,0,0,\ldots,0,0)$,
    $(0,-1,0,\ldots,0,0)$,
    \ldots,
    $(0,0,0,\ldots,-1,0)$.

    Now iteratively construct each polytope $Q_{i+1}$ from the previous
    polytope $Q_i$ as follows:
    \begin{enumerate}[Step 1.]
        \item \label{en-imp-search}
        Search through the visible Euclidean vertices of $Q_i$ for one that
        does not belong to the target polytope $T^\leq$.
        If no such vertex exists (indicating that $Q_i = T^\leq$),
        then output all visible Euclidean vertices
        of $Q_i$ and terminate the algorithm.

        \item \label{en-imp-bdry}
        Otherwise let $\mathbf{v}$ be such a vertex,
        and compute the unique point $\mathbf{x} \in \oproj{p}$ at which
        the line segment joining $\mathbf{v}$ with $\overline{\mathbf{y}}$
        crosses the boundary of the target polytope $T^\leq$.

        \item \label{en-imp-half}
        Compute a projective half-space $H \subseteq \oproj{p}$ that contains
        the target polytope $T^\leq$, and that has the point
        $\mathbf{x}$ in its bounding hyperplane.

        \item \label{en-imp-hull}
        The next polytope is now $Q_{i+1} = Q_i \cap H$.
        Compute a half-space representation and a convex hull
        representation of $Q_{i+1}$.
    \end{enumerate}
\end{algorithm}

\begin{figure}[htb]
    \centering
    \includegraphics[scale=0.3]{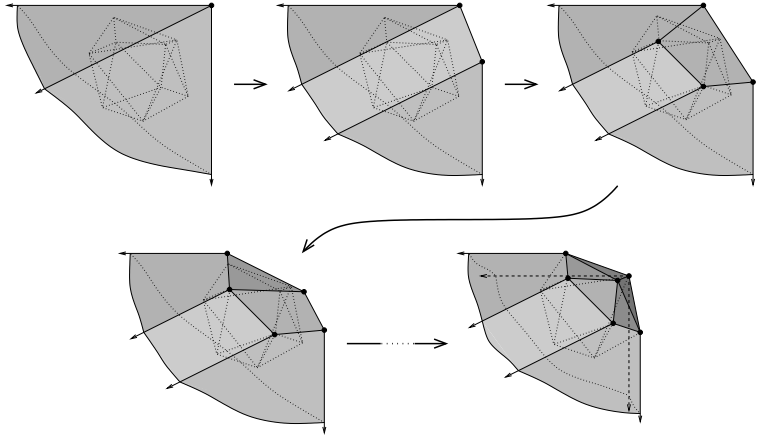}
    \caption{The new algorithm as seen from visible Euclidean space}
    \label{fig-outer-proj}
\end{figure}

The overall procedure is illustrated in Figure~\ref{fig-outer-proj}
(which only shows the portions of $Q_0,Q_1,\ldots$ in visible Euclidean
3-space, and does not show the additional points at infinity).
It should be noted that all of the geometric entities used in this
algorithm lie in the interior of a projective half-space,\footnote{%
    One such half-space is
    $\left\{ (x_1,\ldots,x_{p+1}) \in \oproj{p}\,\left|\,
    x_1 + \ldots + x_p \leq (y_1^\mathrm{max}+\ldots+y_p^\mathrm{max} + 1)
    \cdot x_{p+1} \right.\right\}.$}
which means that all of the relevant concepts (such as polytopes,
convex hulls and line segments) are properly defined.

We will shortly discuss the implementation of the various steps of
this algorithm.  However, our most immediate task is to prove the
correctness of the algorithm as a whole.

\begin{theorem}
    Algorithm~\ref{a-imp} is correct (i.e., it outputs precisely the
    efficient extreme outcomes for our multiobjective linear
    program), and it terminates in finitely many iterations.
\end{theorem}

\begin{proof}
    We prove the correctness of this algorithm in several stages.
    \begin{enumerate}[(i)]
        \item \emph{For each polytope $Q_i$ we have $Q_i \supseteq T^\leq$:}

        This is simple to prove by induction.  The visible Euclidean
        points of the initial simplex $Q_0$ are precisely those with
        signed homogeneous coordinates of the form
        $(y_1^\mathrm{max}-\beta_1,\ldots,y_p^\mathrm{max}-\beta_p,1)$
        for $\beta_1,\ldots,\beta_p \geq 0$, which includes all
        visible Euclidean points of $T^\leq$.  Moreover, the vertices at
        infinity for $Q_0$ are precisely the vertices at infinity for
        $T^\leq$ (Lemma~\ref{l-leq-vertices}).  Thus $Q_0 \supseteq T^\leq$.

        For the inductive step, we observe that $Q_{i+1}$ is the
        intersection of $Q_i$ with a half-space $H \supseteq T^\leq$
        (Step~\ref{en-imp-half}).  It follows that if $Q_i \supseteq T^\leq$
        then $Q_{i+1} \supseteq T^\leq$ also.

        \item \label{en-imp-pf-inf}
        \emph{For each polytope $Q_i$, the vertices at infinity are
        precisely the points
        $(-1,0,0,\ldots,0,0)$,
        $(0,-1,0,\ldots,0,0)$,
        \ldots,
        $(0,0,0,\ldots,-1,0)$:}

        Let $H_\infty$ denote the infinite hyperplane in $\oproj{p}$
        (as defined earlier in Section~\ref{s-proj}).
        The intersection of the initial simplex $Q_0$ with $H_\infty$
        is precisely the $(p-1)$-simplex whose vertices are
        $(-1,0,0,\ldots,0,0)$,
        $(0,-1,0,\ldots,0,0)$,
        \ldots,
        $(0,0,0,\ldots,-1,0)$.
        However, Definition~\ref{d-leq-proj} shows that the intersection
        of $T^\leq$ with $H_\infty$ is this same
        $(p-1)$-simplex.  Because $Q_0 \supseteq Q_1 \supseteq \ldots
        \supseteq T^\leq$, it follows that
        $Q_i \cap H_\infty = T^\leq \cap H_\infty$ for each $i$,
        and in particular the vertices of $Q_i$ at infinity are the same
        $p$ points described above.

        \item \label{en-imp-pf-term}
        \emph{The termination condition is correct, i.e.,
        we terminate if and only if $Q_i = T^\leq$:}

        Step~\ref{en-imp-search} instructs us to terminate if and only
        if no visible Euclidean vertex of $Q_i$ lies outside $T^\leq$.
        By fact~(\ref{en-imp-pf-inf}) above, this is equivalent to the
        condition that no vertex of $Q_i$ lies outside $T^\leq$.
        Because $Q_i$ and $T^\leq$ are both projective polytopes
        with $Q_i \supseteq T^\leq$, this is equivalent to the condition
        that $Q_i = T^\leq$.

        \item \emph{The output is correct, i.e., we output precisely
        the set of efficient extreme outcomes for our multiobjective
        linear program:}

        By fact~(\ref{en-imp-pf-term}) above, if we terminate then we
        have $Q_i = T^\leq$, whereupon Step~\ref{en-imp-search} instructs us
        to output all visible Euclidean vertices of $T^\leq$.
        Corollary~\ref{c-leq-complexity} shows that these are precisely the
        efficient extreme outcomes for our multiobjective linear program.
    \end{enumerate}

    Together these facts show that, if the algorithm terminates, it
    produces the correct results.

    Our final task is to prove that the algorithm terminates in finitely
    many steps.  Let $\partial Q_i$ denote the boundary of the
    projective polytope $Q_i$.  Because $Q_i \supseteq T^\leq$, the
    intersection $\partial Q_i \cap T^\leq$ must be a union of faces
    of $T^\leq$.

    We observe that each intersection $\bdry Q_{i+1} \cap T^\leq$
    is strictly larger than $\partial Q_i \cap T^\leq$, since it contains the
    new point $\mathbf{x}$ identified in Step~\ref{en-imp-bdry}.
    In particular, we know that $\mathbf{x} \notin \partial Q_i$,
    since $\mathbf{x}$ lies strictly between the vertex
    $\mathbf{v}$ and the interior point $\overline{\mathbf{y}}$.
    However, we also know that $\mathbf{x} \in \partial Q_{i+1} \cap T^\leq$,
    since Step~\ref{en-imp-bdry} ensures that
    $\mathbf{x} \in T^\leq$ and Step~\ref{en-imp-half}
    ensures that $\mathbf{x} \in \partial Q_{i+1}$.

    It follows that each intersection $\partial Q_{i+1} \cap T^\leq$
    contains strictly more faces of $T^\leq$ than
    the previous intersection $\partial Q_i \cap T^\leq$.
    Because $T^\leq$ has finitely many faces, this sequence cannot grow
    forever and the algorithm must terminate.
\end{proof}

The individual steps of Algorithm~\ref{a-imp} are performed in much the
same way as for the original algorithm.

\begin{itemize}
    \item In the initialisation phase, the point $\mathbf{y}^\mathrm{max}$
    is straightforward to compute, and for the interior point
    $\overline{\mathbf{y}} \in T^\leq$ we can simply use
    $\widehat{\mathbf{y}}$ from Definition~\ref{d-leq-square}.
    The half-space representation of the initial simplex $Q_0$ involves
    the $p$ ``ordinary'' half-spaces with
    signed homogeneous coordinates
    $(-1,0,0,\ldots,0,y_1^\mathrm{max})$,
    $(0,-1,0,\ldots,0,y_2^\mathrm{max})$,
    \ldots,
    $(0,0,0,\ldots,-1,y_p^\mathrm{max})$,
    as well as
    the visible half-space $(0,0,0,\ldots,0,1)$.

    \item In Step~\ref{en-imp-search} we are required to test whether a
    visible Euclidean vertex belongs to $T^\leq$.  For this we can
    restrict our attention to visible Euclidean space (where $T^\leq$
    becomes the Euclidean polyhedron $Y^\leq$), and we can employ
    precisely the same feasibility test used in the original algorithm.

    \item Likewise, Steps~\ref{en-imp-bdry} and~\ref{en-imp-half} can be
    performed in visible Euclidean space using the same feasibility tests,
    univariate search techniques and dual linear programs as before.
    In Step~\ref{en-imp-half} it is trivial to convert the
    resulting half-space $H$ back into oriented projective $p$-space.
\end{itemize}

This leaves Step~\ref{en-imp-hull}, where we must enumerate the
vertices of the new polytope $Q_{i+1}$.
For this we can use the double description method directly
in oriented projective $p$-space:

\begin{theorem} \label{t-dd-proj}
    The results of Theorem~\ref{t-dd} and Lemma~\ref{l-adj} also
    hold true in oriented projective $p$-space.
    That is, they remain true if we replace all occurrences of
    $\euc{p}$, polytopes, half-spaces and hyperplanes with
    $\oproj{p}$, projective polytopes, projective half-spaces and
    projective hyperplanes respectively.
\end{theorem}

\begin{proof}
    Instead of offering a direct proof, we can transform these results
    into their Euclidean counterparts using a powerful tool known as a
    \emph{projective mapping}.  Projective mappings are geometric
    transformations in $\oproj{p}$ that preserve collinearity and concurrency.
    They map projective hyperplanes to projective hyperplanes
    and projective half-spaces to projective half-spaces (and therefore
    projective polytopes to projective polytopes).
    Every projective mapping has an inverse projective mapping, and
    for any projective half-space $H \subseteq \oproj{p}$
    there is a projective mapping that takes $H$ to the visible half-space.
    For more information on projective mappings, the reader is referred
    to Boissonnat and Yvinec \cite{boissonnat98-geometry}.

    Using projective mappings, the proof of Theorem~\ref{t-dd-proj}
    becomes simple.  If $P \subseteq \oproj{p}$ is a projective polytope
    then it lies in the interior of some projective half-space
    $H \subseteq \oproj{p}$.  Let $\pi$ be a projective mapping that
    maps $H$ to the visible half-space.  The interior points of
    $H$ map to visible Euclidean points, and $P$ therefore maps to a
    visible Euclidean polytope.

    We can now apply Theorem~\ref{t-dd} and Lemma~\ref{l-adj} directly
    to the Euclidean polytope $\pi(P)$.  Because both results are
    expressed purely in terms of lines, intersections, half-spaces and
    hyperplanes---all notions that are preserved under a projective
    mapping---the corresponding results are therefore true
    for the original projective polytope $P$.
\end{proof}

When using the double description method in this way to enumerate the
vertices of $Q_{i+1}$, it is important to take into account \emph{all} of
the vertices of $Q_i$, including the vertices at infinity.

This final point highlights a key motivation for working in oriented
projective geometry.
In Euclidean space the double description method works
only for polytopes, not unbounded polyhedra, because it requires a description
of the polytope as a convex hull of its vertices.
Oriented projective geometry allows us to use the double description
method with unbounded Euclidean polyhedra (such as $Y^\leq$) by treating
them as polytopes---that is, convex hulls of vertices---in $\oproj{p}$.

\subsection{Running time} \label{s-outer-complexity}

The primary benefit of our new outer approximation algorithm is
that we avoid constructing and analysing a large number of ``unnecessary''
vertices that are \emph{not} efficient extreme outcomes for our
multiobjective linear program.  As shown in Section~\ref{s-eff}, the
projective polytope $T^\leq$ (used in the new algorithm) has only
$p$ unnecessary vertices, which are the $p$ vertices at infinity.  In
contrast, the Euclidean polytope $Y^\square$ (used in the original
algorithm) has a significant number of unnecessary vertices---their
number is always exponential in the number of objectives $p$
(Corollary~\ref{c-square-best}), and can also grow at a severe rate
relative to the problem size $(n,m)$ (Corollary~\ref{c-square-worst}).

Just as $Y^\square$ has significantly more vertices than $T^\leq$,
we also expect the intermediate polytopes $P_1,P_2,\ldots$ from the
original algorithm to have significantly more vertices than
$Q_1,Q_2,\ldots$ from the new algorithm (recall that these
families of polytopes successively approximate $Y^\square$ and $T^\leq$
respectively).
This growth in the number of vertices has a direct impact on the running
time of the outer approximation algorithm (Algorithm~\ref{a-outer}):
\begin{itemize}
    \item In Step~\ref{en-outer-search}, the time required to search for
    a vertex $\mathbf{v}$ outside the target polytope $Y^\square$ is directly
    proportional to the number of vertices in the intermediate polytope $P_i$.

    \item In Step~\ref{en-outer-hull}, the time required to run the
    double description method is proportional to the \emph{cube} of the
    number of vertices:  if the intermediate polytope $P_i$ has $V$ vertices,
    we have $O(V^2)$ potential pairs
    $\mathbf{u} \in S_+$ and $\mathbf{v} \in S_-$
    (Theorem~\ref{t-dd}), and for each such pair the adjacency test
    requires us to search through all of the vertices again
    (Lemma~\ref{l-adj}).

    \item More generally, the number of stages in the algorithm as a whole
    (that is, the number of intermediate polytopes $P_1,\ldots,P_k$) depends
    on the number of vertices:  each stage corresponds to some intermediate
    vertex $\mathbf{v}$ outside the target polytope $Y^\square$,
    and so with more intermediate vertices we expect more intermediate
    polytopes in total.
\end{itemize}

All of these observations indicate that the new algorithm
should offer significant benefits over the original algorithm
in terms of running time.

The only added complication for the new algorithm
is that we must perform arithmetic in oriented projective $p$-space,
where we work with $p+1$ coordinates at a time instead of $p$.
The cost of this ``extra coordinate'' is negligible, and (unlike
the number of vertices) does not affect the
asymptotic growth of the running time at all.

%% file: conc.tex
\section{Conclusions} \label{s-conc}

In this paper we show how oriented projective geometry can be
applied to the field of multiobjective linear programming, allowing us
to work with efficiency-equivalent \emph{projective} polytopes
that are significantly less complex than their Euclidean counterparts.
Polytopes (as opposed to unbounded polyhedra) are important because they
support techniques such as the double description method, which
operates on convex hulls of finitely many vertices.
The reduction in complexity is important because it leads to a
significant reduction in running time.

As a concrete illustration, we apply these techniques to the outer
approximation algorithm, which is used to generate all efficient
extreme outcomes for a multiobjective linear program.  By working in
oriented projective geometry we obtain a considerable reduction in
the number of \emph{non-efficient} vertices that the algorithm
generates, and we show how this directly benefits the
running time as a result.

We support these results with explicit asymptotic estimates on the
number of non-efficient vertices generated in the old Euclidean
setting and the new oriented projective setting.  Using oriented projective
geometry, the number of non-efficient vertices is precisely the number
of objectives $p$.  In the old Euclidean setting this number is at least
$2^p-1$, and it can grow as fast as $\Omega(n^{\lfloor (p-1)/2 \rfloor})$.

For applications with even a moderate number of objectives, this growth
rate can be extremely high.  For example, Steuer
describes a blending problem with $p=8$ objectives
\cite{steuer84-sausage}; here the
number of non-efficient vertices in the old algorithm
can grow like $\Omega(n^3)$, which for
large problem sizes is a significant burden.
When the number of objectives is higher, as in the
aircraft control design problem of Schy and Giesy with $p=70$ objectives
\cite{schy88-aircraft},
even the best-case growth rate of $2^p-1$ becomes crippling.
In contrast, the new techniques in this paper deliver just $8$ and $70$
non-efficient vertices respectively, which is negligible in comparison.

Although we focus on the outer approximation algorithm for our
application in Section~\ref{s-outer}, the underlying techniques of
Sections~\ref{s-proj} and~\ref{s-eff} are quite general.
We expect that the methods described in this paper can find many
fruitful applications throughout the field of multiobjective linear
programming, and that with further research the geometric and analytical
insights described here can lead to new discoveries in other areas
of high-dimensional optimisation.